\documentclass[twocolumn]{IEEEtran}  

\usepackage{graphicx}      

\usepackage[table,xcdraw]{xcolor}
\usepackage{amsmath}

\usepackage{algorithm, algorithmic}
\usepackage{amsthm,amsmath,amssymb,mathrsfs}

\usepackage{tasks}
\usepackage{bbold}
\usepackage[utf8]{inputenc}
\usepackage{dsfont}
\usepackage{adjustbox}
\usepackage{mathtools, cuted}
\usepackage{graphics}
\usepackage{support-caption}
\usepackage{subcaption}
\usepackage{enumitem}
\usepackage{support-caption}
\usepackage{subcaption}

\usepackage{multirow}

\usepackage{booktabs}
\usepackage{threeparttable}

\newtheorem{theorem}{\emph{Theorem}}
\newtheorem{proposition}{Proposition}
\newtheorem{lemma}{Lemma}
\newtheorem{remark}{Remark}

\newtheorem{definition}{Definition}

\title{\LARGE \bf
\textcolor{black}{\color{red}Robust optimal policies for team Markov games}}
\author{Feng Huang, Ming Cao, and Long Wang
\thanks{This work was supported by the National Natural Science
Foundation of China (grant no. 62036002). F.H. acknowledges the support from China Scholarship Council (grant no. 201906010075). M.C. was supported in part by the European
Research Council (grant no. ERC-CoG-771687). (Corresponding authors: Long Wang and Ming Cao)}
\thanks{ F. Huang is with the Center for Systems and Control, College of Engineering, Peking University, Beijing 100871, China and
        ENTEG, Faculty of Science and Engineering,
University of Groningen, Groningen 9747 AG, The Netherlands (e-mail: fenghuang@pku.edu.cn).
    }
\thanks{ M. Cao is with
        ENTEG, Faculty of Science and Engineering,
        University of Groningen, The Netherlands (e-mail: m.cao@rug.nl).
    }
\thanks{
L. Wang is with the Center for Systems and Control, College of Engineering, Peking University, Beijing 100871, China (e-mail: longwang@pku.edu.cn).
}
}

\begin{document}
\maketitle
 \thispagestyle{plain}
\pagenumbering{arabic}
 \pagestyle{plain}

\begin{abstract}
In stochastic dynamic environments, team {\color{red}Markov} games have emerged as a versatile paradigm for studying sequential decision-making problems of fully cooperative multi-agent systems. However, the optimality of the derived policies is usually sensitive to model parameters, which are typically unknown and required to be estimated from noisy data in practice. To mitigate the sensitivity of optimal policies to these uncertain parameters, we propose a ``robust'' model of team {\color{red}Markov} games in this paper, where agents utilize \emph{robust optimization} approaches to update strategies. This model extends team {\color{red}Markov} games to the scenario of incomplete information and meanwhile provides an alternative solution concept of robust team optimality. To seek such a solution, we develop {\color{red}a robust iterative learning algorithm of team policies} and prove its convergence. This algorithm, compared with robust dynamic programming, not only possesses a faster convergence rate, but also allows for using approximation calculations to alleviate the curse of dimensionality. Moreover, some numerical simulations are presented to demonstrate the effectiveness of the algorithm by generalizing the game model of sequential social dilemmas to uncertain scenarios.

\end{abstract}

\begin{IEEEkeywords}
Robust optimal control, team Markov games, multi-agent learning, sequential social dilemmas.
\end{IEEEkeywords}


\section{INTRODUCTION}
\IEEEPARstart{M}{arkov} games~\cite{littman1994markov}, also known as stochastic games~\cite{shapley1953stochastic}, as a general framework of multi-agent reinforcement learning (MARL) have long been an active research topic across the fields of stochastic optimal control, operations research, and artificial intelligence (AI)~\cite{busoniu2008comprehensive,bertsekas2019reinforcement,zhang2019multi}. Especially, it has attracted much attention in recent years due to some ground-breaking advances made by the MARL in conjunction with deep neural networks in achieving human-level performance in several long-standing sequential decision-making conundrums~\cite{jaderberg2019human,vinyals2019grandmaster}. Compared with the normal-form matrix games~\cite{osborne1994course,bacsar1998dynamic} and evolutionary games~\cite{smith1982evolution,riehl2018survey}, one of its most distinctive features is
the introduction of {\color{red}game states.}
In each stage of the play, the game is in a {\color{red}state} taken from a given set, and each player then relies on its decision rule as a function of the current {\color{red}state} to choose actions. The collection of players' actions, together with the current {\color{red}state}, subsequently determine the probability distribution of the state that the play will visit at next stage{\color{red}. As} the consequence of the joint actions and {\color{red}state} transitions, every player will receive an immediate payoff. Therefore, Markov games are usually deemed as the generalization of one-shot matrix games to the dynamic multi-stage and multi-state situations, in which {\color{red}game states} change in response to players' decisions. {\color{red} Also, they extend the Markov decision processes (MDPs)~\cite{puterman2014markov} from only involving a single agent to the scenario involving multiple decision-makers~\cite{solan2015stochastic}.}
\par
In contrast to the non-cooperative setting of Markov games, team Markov games~\cite{wang2002reinforcement,matignon2012} are a fully cooperative multi-agent system, in which a group of players works together, through interactions, coordinations, and information-sharing, to collectively accomplish a task or achieve a goal~\cite{boutilier1996planning,panait2005cooperative}. Hence, it is formally defined as the Markov games where all players have a common payoff function~\cite{busoniu2008comprehensive,zhang2019multi}.
In particular, such a game model has recently sparked
increasing research interest across various disciplines due to its extensive real-world applications, such as multi-robot systems, unmanned aerial vehicles, and communication networks (see~\cite{busoniu2008comprehensive,zhang2019multi,oroojlooyjadid2021review} for an overview), and its close connection to the theory of team decisions~\cite{ho1972team,marden2018game,bora21}. However, most of the existing work is based on the complete information setting where the structural information of the game, such as players' payoff functions and state transition probabilities, is the common prior knowledge to all players. Under such an assumption, the notion of the optimal Nash equilibrium~\cite{wang2002reinforcement}, i.e. the joint decision rule achieving the maximal expected discounted sum of the gains received by the team, naturally provides a plausible solution to predicting the outcome of the game evolution. This is because at this equilibrium state, no player can improve the long-term expected return of the team by unilaterally deviating from its policy. Although there are some debates about the predictive or prescriptive role of the equilibrium concept in practice (see~\cite{shoham2007if} and its commentaries), Nash equilibria have played a central role in the corresponding studies~\cite{salehisadaghiani2016distributed,lu2020online}.
\par
While the Nash equilibrium solution always exists for the Markov games with complete information~\cite{fink1964equilibrium}, in the real world, the problem related to the model uncertainty of games is pervasive. For example, in many realistic applications associated with reinforcement learning (RL), the reward functions of agents are usually required to be designed or learned from interactions, and their properties strongly affect the success of targeted tasks~\cite{dulac2021challenges}; also, the transition probability distribution of states is generally estimated from historical data, thereby influenced by statistical errors~\cite{mannor2019data}. In particular, such an issue has given rise to well-grounded concerns, such as in the field of AI safety~\cite{amodei2016concrete,huang2020survey} and uncertain robotic systems~\cite{fisac2018general}, and accordingly it prompts the research priorities of robust AI~\cite{russell2015research}.
\par
\subsection{\color{red} Related work}
{\color{red}In} game theory, {\color{red}the problem related to model uncertainty} has a long research tradition and {\color{red}usually refers to} the games with ``incomplete information"~\cite{harsanyi1967games}, {\color{red}
in which part of game parameters
is unspecified.
}
The first analytical framework developed to study this class of games is from the seminal three-part essay by Harsanyi~\cite{harsanyi1967games} where a new model named ``Bayesian games'' is constructed and the notion of Nash equilibrium is extended to the incomplete information scenario, termed ``Bayesian equilibrium.''
Utilizing robust optimization~\cite{ben2009robust},
Aghassi and Bertsimas then relax the prior distribution assumption of Harsanyi's model
for uncertain parameters and open a new avenue to robust games~\cite{aghassi2006robust}, where players take a best response depending on a worst-case payoff matrix in one-shot normal-form games. Via a similar technique, these results are later extended to Markov games~\cite{kardecs2011discounted} and some deep MARL algorithms have been developed to find the equilibria~\cite{zhang2020robust}.
\par
{\color{red}Parallel} to the game-theoretical research, robust dynamic programming (rDP)~\cite{nilim2005robust,iyengar2005robust} is another domain devoting to coping with the decision-making problems with data uncertainty. To mitigate the sensitivity of optimal policies to the ambiguity of transition probabilities in MDPs, two algorithms, the so-called robust value iteration (rVI) and robust policy iteration (rPI), are proposed to find robust optimal policies. Inspired by these work, Kaufman and Schaefer~\cite{kaufman2013robust} then develop a robust version of modified policy iteration (rMPI)~\cite{puterman2014markov} as the generalization of these two algorithms.  {\color{red}In particular, recent years have witnessed some applications of these methods in the temporal
logic control of stochastic systems~\cite{sofie21} and the continuous control of physical systems~\cite{lutter2021robust}.}
However, one common shortcoming of these algorithms is that the convergence rate may be quite slow when the discount factor used to compute the long-term expected return is close to one.
{\color{red}As the policy space is very large and/or continuous, moreover, these algorithms may become less applicable in general due to the curse of dimensionality, and thus require to resort to function approximations or data-driven techniques~\cite{mannor2019data}. Motivated by this fact, a large body of robust RL algorithms has recently been proposed by incorporating diverse approaches, such as adversarial training~\cite{pinto2017robust,phan2021resilient}, online policy search~\cite{Mankowitz2020Robust,wang2021online}, and least squares policy iteration~\cite{badrinath21a}.
}
\subsection{\color{red} Contributions}
In this work, we aim to deal with {\color{red}the sequential decision-making problem of a cooperative team in stochastic uncertain environments,} and accordingly propose a robust model of team {\color{red} Markov} games. This model relaxes the complete information assumption in team {\color{red} Markov} games and meanwhile provides an alternative solution concept termed robust
team-optimal policy. To seek such a solution, we develop {\color{red}a robust iterative learning algorithm of team policies, which we call \emph{robust approximate Team Policy Iteration} (raTPI).} Compared with rDP, this algorithm allows for using approximation computations to alleviate the curse of dimensionality. Furthermore, we present the convergence analysis of the algorithm within mild approximation tolerance and accordingly calculate its convergence rates. The results manifest that {\color{red}at a near exponential convergence rate}, this algorithm can effectively find the robust team-optimal policy {\color{red}with sufficient accuracy} after a finite number of iterations. To demonstrate the effectiveness of the algorithm, {\color{red}we also carry out numerical simulations by} generalizing the canonical game model of social dilemmas to sequential uncertain scenarios. Using them as benchmarks, simulations show that, as compared with rVI and rMPI, our algorithm admits substantially less iteration time to find the robust team-optimal policy.
\par
\textbf{Notation}: Throughout the paper, we use $\mathbb{R}$ and $\mathbb{N}$ to represent the set of real numbers and the set of non-negative integers, respectively. For $n$ sets $\mathcal{A}^i$, $i=1,2,\ldots,n$, let $\ltimes_{i=1}^{n}\mathcal{A}^i$ be their Cartesian product. For a scalar $x$, we use $|x|$ to denote its absolute value. The space of {\color{red}probability distributions} on the set $\mathcal{S}=\{s^1,s^2,...,s^m\}$ is denoted by $\Delta(\mathcal{S})$, and the set of bounded real valued functions on $\mathcal{S}$ is represented by $\mathcal{V}$. For a vector $v=[v(s^1),v(s^2),\ldots,v(s^m))]^T\in \mathcal{V}$, its vector norm is defined by $\|v\|:=\sup_{s\in \mathcal{S}}|v(s)|$, where $v(s)$ represents the component of $v$ corresponding to $s\in \mathcal{S}$ and $T$ represents transpose. As such, $(\mathcal{V},\|\cdot\|)$ forms a normed linear space, and it is also a Banach space. For a matrix $P=[p_{kl}]\in \mathbb{R}^{m\times m}$ with element $p_{kl}$ in row $k$ and column $l$, its matrix norm and spectral radius are defined by $\|P\|:=\sup_{k\in\{1,2,\ldots,m\}}
\sum_{l=1}^{m}|p_{kl}|$ and $\sigma(P):= \limsup_{t\rightarrow \infty}
\|P^t \|^{1/t}$, respectively. Alternatively, we also use $P(l\mid k)$ to represent the element of $P$ in row $k$ and column $l$. We use $I$ and $\mathbf{1}$ to denote the identity matrix and the all-ones vector with appropriate dimensions, respectively. {\color{red}The order, maximum, and minimum of vectors and matrices refer to componentwise order, maximum, and minimum}.

\section{Problem formulation and preliminaries}
{\color{red}In this section, we first introduce the sequential decision-making problem in team Markov games and then present one of its robust counterparts in uncertain environments.
}
\subsection{\color{red}Team Markov games}
{\color{red}We consider the decision-making problem of a cooperative team with $n$ players, which is modelled by a team Markov game~\cite{wang2002reinforcement} (also known as multi-agent MDPs~\cite{boutilier1996planning}) with the infinite decision horizon $\mathcal{T}:=\{1,2,\ldots\}$. A team Markov game is characterized by a tuple $$\langle\mathcal{N},\mathcal{S},\mathcal{A},
\{r^i\}_{i\in\mathcal{N}},p\rangle,$$
where $\mathcal{N}:=\{1,2,\ldots,n\}$ is the finite set of the indices of $n$ players; $\mathcal{S}:=\{s^1,s^2,...,s^m\}$ is the finite set consisting of $m$ states shared by all players; $\mathcal{A}:=\ltimes_{i=1}^{n}\mathcal{A}^i$ represents the set of the joint action of all players, which is the Cartesian product of the discrete finite set $\mathcal{A}^i$ of actions available to player $i\in \mathcal{N}$; $r^i(s,a,s'):\mathcal{S}\times\mathcal{A}\times\mathcal{S}\rightarrow\mathbb{R}$
defines the payoff function of player $i$, which
depends on the current state $s$, the joint action $a$ of all players, and the state $s'$ at the next epoch; and $p(\cdot|s,a):\mathcal{S}\times \mathcal{A}\rightarrow \Delta(\mathcal{S})$ represents the transition probability function of states, which maps from the current state $s$ to the probability distribution over the state space, given the joint action $a$.
}
\par
Without loss of generality, we assume in this paper that both the payoff functions of players and the transition probability functions of {\color{red}states} are stationary. That is, for any triple $(s,a,s')\in \mathcal{S}\times\mathcal{A}\times\mathcal{S}$, $r^i(s,a,s')$ for $\forall i \in \mathcal{N}$ and $p(s'|s,a)$ do not change with time. Moreover, $r^i(s,a,s')$ is assumed to be bounded, i.e. $|r^i(s,a,s')|\leq R_{\max}<\infty$ for $\forall i\in \mathcal{N}$ and $\forall (s,a,s')\in \mathcal{S}\times\mathcal{A}\times\mathcal{S}$, where $R_{\max}$ is a constant.
\par
{\color{red}
Moreover, in this paper, we restrict our attention to stationary and deterministic policies due to their mathematical tractability~\cite{puterman2014markov,bertsekas1996neuro}. Specifically, at each decision epoch, we consider that given the current state $s\in \mathcal{S}$, each player $i\in\mathcal{N}$ in the team selects an action $a^i\in \mathcal{A}^i$ deterministically according to its individual decision rule $d^i(s):\mathcal{S} \rightarrow \mathcal{A}^i$, which is a mapping from state $s$ to a deterministic action $a^i$. These decision rules adopted by
player $i$ at all decision epochs then constitute its policy or strategy used in the game, which is a time sequence $(d^i,d^i,\ldots)$ of decision rules. For ease of notations, we denote the policy of player $i$ by $\pi^i:=(d^i)^\infty=(d^i,d^i,\ldots)$. As such, the team decision rule and the team policy can be represented by
the joint decision rule of all players, $d:=(d^1,\ldots,d^n)$, and the joint policy of all players, $\pi:=(\pi^1,\ldots,\pi^n)=(d)^\infty$, respectively. Accordingly, we denote the space of team decision rules and the space of team policies by $\mathcal{D}$ and $\Pi$, respectively. From the definition of team decision rules, one can see that every joint action of all players is completely determined by the team decision rule $d$.
In response to the joint action $a_\tau$ given by $d$ at decision epoch $\tau\in\mathcal{T}$, the game system will transit from the current state $s_\tau$ to the state $s_{\tau+1}$ at the next epoch with probability $p(s_{\tau+1}|s_\tau,a_\tau)$. As the consequence of the joint action and the state transition, every player $i\in\mathcal{N}$ will then get an immediate payoff $r^i(s_\tau,a_\tau,s_{\tau+1})$.
\par
Starting from an initial state $s_1=s\in \mathcal{S}$, when the team adopts a specific joint policy $\pi\in \Pi$ to play the game, the system will induce a probability measure on the trajectory of states and the sequence of joint actions. To evaluate the performance of the joint policy, we adopt the expected total discounted return of the team-average payoff as the criterion, which is calculated by
\begin{equation}\label{eq8}
 v_{\pi}(s)=\mathbb{E}_{\pi}\left\{ \sum_{\tau=1}^{\infty}\lambda^{\tau-1}
  r(s_\tau,a_\tau,s_{\tau+1})\bigg|s_1=s \right\},
\end{equation}
where $\mathbb{E}_{\pi}\{\cdot\}$ represents the expectation over the stochastic process $\{(s_\tau,a_\tau)\}_{\tau\in \mathcal{T}}$ induced by the joint policy $\pi$; $\lambda\in [0,1)$ is the discount factor used to evaluate the present value of future payoffs; and $r(s_{\tau},a_{\tau},s_{\tau+1})
=\frac{1}{n}\sum_{i=1}^{n}r^i(s_{\tau},a_{\tau},s_{\tau+1})$ is the team-average payoff at epoch $\tau\in \mathcal{T}$. As a solution to the team decision-making problem, we consider that all players aim to seek a \emph{team-optimal policy} $\pi'\in \Pi$ such that for $\forall \pi\in \Pi$,
$$
v_{\pi'}(s)\geq v_{\pi}(s), \forall s \in \mathcal{S}.
$$
\par
To find $\pi'$, one can see from~(\ref{eq8}) that all players will rely on the common utility $r$ to make decisions at each epoch. In other words, all players in effect have a common payoff function $r$ in the game. Thus, the team decision-making problem is in a fully cooperative setting~\cite{busoniu2008comprehensive,zhang2019multi}. Moreover, it is easy to verify that every team-optimal policy is a Markov perfect equilibrium (see \cite{zhang2019multi,kardecs2011discounted} for its definition) of Markov games.

}

For ease of the exposition of our main results, we further introduce some vector and matrix notations, and rewrite~(\ref{eq8}) in a vector form. Given a joint action $a\in \mathcal{A}$, let $P_{a}$ be the transition probability matrix of {\color{red}states}, where {\color{red}the element of $P_a$ in row $k$ and column $l$ is given by $p(s^l|s^k,a), \forall s^k,s^l\in\mathcal{S}$.} Accordingly, we denote the transition probability matrix induced by the joint decision rule $d\in \mathcal{D}$ by $P_{d}$, where the entry of $P_{d}$ in row $k$ and column $l$ is given by $p(s^l|s^k,d(s^k))$.
As such, given a joint policy $\pi=(d)^\infty$, the probability $Pr(s_{\tau}=s^l|s_{1}=s^k;\pi)$ that the game system transits from the initial state $s^k\in\mathcal{S}$ to the state $s^l\in\mathcal{S}$ at epoch $\tau\in \mathcal{T}$ can be calculated by
\begin{equation}\label{eq10}
  Pr(s_{\tau}=s^l|s_{1}=s^k;\pi)
  =\left(P_{d}P_{d}
  \cdots P_{d}\right)(l|k)=
  \left(P_{d}\right)^{\tau-1}(l|k),
\end{equation}
where $(P_{d})^{0}=I$. {\color{red}For a given $P_d$, moreover, one can} calculate the expected value $r_{(d,P_d)}(s)$ of the team-average payoff for $\forall s\in\mathcal{S}$ by
\begin{equation}\label{eq11}
    r_{(d,P_d)}(s)=
     \sum_{s'\in\mathcal{S}}r(s,d(s),s')p(s'|s,d(s)).
\end{equation}
{\color{red}It follows from~(\ref{eq10}) and (\ref{eq11}) that for $s_1=s^k\in \mathcal{S}$, (\ref{eq8}) can be rewritten by}
\begin{equation}\label{eq4}
\begin{split}
v_{\pi}(s^k)&=\sum_{\tau=1}^{\infty}\lambda^{\tau-1}
\mathbb{E}_{\pi}\{r(s_\tau,a_\tau,s_{\tau+1})|s_1=s^k\} \\
&=\sum_{\tau=1}^{\infty}\lambda^{\tau-1}
\sum_{l=1}^{m}
Pr(s_{\tau}=s^l|s_{1}=s^k;\pi)r_{(d,P_d)}(s^l) \\
&=\sum_{\tau=1}^{\infty}\lambda^{\tau-1}
\sum_{l=1}^{m}
\left(P_{d}\right)^{\tau-1}(l|k)r_{(d,P_d)}(s^l).
\end{split}
\end{equation}
\par
Let $v_{(d,P_d)}:=[v_{\pi}(s^1),v_{\pi}(s^2),\ldots,v_{\pi}(s^m)]^T$ and $r_{(d,P_d)}:=[r_{(d,P_d)}(s^1),r_{(d,P_d)}(s^2),
\dots, r_{(d,P_d)}(s^m)]^T$.
Then, from~(\ref{eq4}), the vector form of~(\ref{eq8}) can be given by
\begin{equation}\label{eq5}
\begin{split}
   v_{(d,P_d)}&=\sum_{\tau=1}^{\infty}\lambda^{\tau-1}
   {P_{d}}^{\tau-1}r_{(d,P_d)} \\
   &=r_{(d,P_d)}+\lambda P_{d} \left(r_{(d,P_d)}+\lambda P_{d}r_{(d,P_d)}+\cdots\right)\\&
   =r_{(d,P_d)}+\lambda P_{d}v_{(d,P_d)}.
\end{split}
\end{equation}
Since $\|r_{(d,P_d)}\|\leq R_{\max}$ in view of the assumption
$|r^i(s,a,s')|\leq R_{\max}<\infty$ for $\forall i\in \mathcal{N}$ and $\forall (s,a,s')\in \mathcal{S}\times\mathcal{A}\times\mathcal{S}$ and $P_d$ is a row stochastic matrix, $v_{\pi}(s)\leq R_{\max}/(1-\lambda)$ for $\forall s\in \mathcal{S}$ and $\forall \pi \in \Pi$. It implies that $v_{(d,P_d)}\in \mathcal{V}$ for $\forall d\in \mathcal{D}$ and any $P_d$. Moreover, note that $\sigma(\lambda P_{d})\leq \|\lambda P_{d}\|=\lambda<1$. It then follows that $v_{(d,P_d)}=(I-\lambda P_{d})^{-1}r_{(d,P_d)}
=\sum_{\tau=1}^{\infty}\left(\lambda P_{d}\right)^{\tau-1} r_{(d,P_d)}$ is the unique solution to the equation
$
v=r_{(d,P_d)}+\lambda P_{d}v,
$
where $v\in \mathcal{V}$ is the variable. For the sake of convenience, we will call any vector $v\in \mathcal{V}$ the value function afterwards.

\subsection{Robust team Markov games}
The aforementioned model of team Markov games is based on the assumption of complete information, i.e. the parameter information of the game is explicitly known to players. Here, we propose one of its robust counterparts by relaxing this assumption. We assume that players do not know the true transition probabilities of game {\color{red} states}, but rather commonly perceive an uncertainty set of their possible values. More specifically, given a joint action $a\in \mathcal{A}$, the transition probability matrix $P_{a}$ is not predetermined but lies in an uncertainty set $\mathcal{P}_{a} \subset \mathbb{R}^{m\times m}$. Since players' payoff functions are defined to depend on the state at the next epoch in team Markov games,
the fluctuation of transition probabilities will lead to the uncertainty of
 players' payoffs. Therefore, we refer to a team Markov game as a \emph{robust team Markov game} if players' payoffs and/or the transition probabilities of {\color{red}states} are uncertain.
\par
Following the convention, we assume that for $\forall a\in \mathcal{A}$, the uncertainty set $\mathcal{P}_{a}$ satisfies the so-called ``rectangularity" property~\cite{nilim2005robust,iyengar2005robust}. Namely, $\mathcal{P}_{a}$ has the form of $\mathcal{P}_{a}=\ltimes_{k=1}^{m}\mathcal{P}_{a}(\cdot|k)$ and is independent of historically visited states and actions, where $\mathcal{P}_{a}(\cdot|k)\subseteq \Delta(\mathcal{S})$ for $\forall k\in\{1,2,\ldots,m\}$ characterizes the uncertainty of the $k$-th row of $P_a$ {\color{red} and is assumed to be discrete and finite for ease of exposition. The extension to the continuous case is straightforward, which will not significantly bring new insights but rather require more complex notations and additional assumptions for the existence of optimal policies~\cite{puterman2014markov,kaufman2013robust}. Moreover, in many practical applications, such as video games, board games, and autonomous vehicles, a robust learning system acquires all possible transition probabilities of different scenarios usually by massively sampling from the realistic physical system, even though the underlying uncertainty set may be continuous. Therefore, in practice, the learning system generally operates in the discrete and finite sample set. Also, such a discrete and finite set is more practical for the numerical evaluation of an algorithm.} Under this assumption, the admissible set of $P_d$ for a given $d\in\mathcal{D}$ can be then represented by $\mathcal{P}_d:=\ltimes_{k=1}^m\mathcal{P}_{d(s^k)}(\cdot|k)$. To systematically mitigate the influence of uncertain transition probabilities on the performance of the team-optimal policy, we consider the \emph{robust team-optimal policy} $\pi^*=(d^{*})^{\infty}\in \Pi$ as an alternative solution such that the value function $v_{(d^*,P_{d^{*}}^*)}$ is maximal with respect to the worst-case $P_{d^{*}}^*$ in the uncertainty set $\mathcal{P}_{d^*}$.
\begin{definition}
A value function $v^*:=v_{(d^*,P_{d^*}^{*})}$ is said to be robust team-optimal if $v_{(d^*,P_{d^*}^{*})}=\max_{d\in \mathcal{D}} \min_{P_d\in \mathcal{P}_d} v_{(d,P_d)}$, and accordingly $\pi^*=(d^{*})^{\infty}$ is called a robust team-optimal policy.
\end{definition}
Similar to the team-optimal policy, it is also easy to verify that every robust team-optimal policy is a robust Markov perfect equilibrium (see~\cite{kardecs2011discounted} for its definition) of robust stochastic games. Note that both $\mathcal{D}$ and $\mathcal{P}_d$ for $\forall d\in \mathcal{D}$ are discrete and finite. Therefore, such a robust team-optimal policy can always be attained. Moreover, in contrast to the robust team-optimal policy, its near-optimal counterpart may be more preferred in practice.
{\color{red}
\begin{definition}
A policy $(d^{\epsilon})^\infty$ is said to be \emph{$\epsilon$-robust team-optimal} for $\epsilon>0$ if
  $\min_{P_{d^{\epsilon}}\in \mathcal{P}_{d^{\epsilon}}} v_{(d^{\epsilon},P_{d^{\epsilon}})}\geq v^*-\epsilon\mathbf{1}$.
\end{definition}
}
\par
In the following section, {\color{red}we will develop an algorithm to seek the robust team-optimal policy and meanwhile present its convergence analysis.}

\section{\color{red}Robust team-optimal policy seeking}
To seek the robust team-optimal policy, we here develop {\color{red}an algorithm, named robust approximate Team Policy Iteration (raTPI). This algorithm is built on the structure of the so-called ``classical
information pattern''~\cite{witsenhausen1968counterexample,bertsekas2021multiagent} or the framework of centralized learning with decentralized execution, under which a central controller is used to compute the desired team policy $\pi$ and then each individual policy of $\pi$ is communicated to the corresponding player for execution.
}
\subsection{The raTPI Algorithm}

{\color{red}
The basic idea of the raTPI algorithm is to utilize rMPI to unify the advantages of both rVI and rPI, and meanwhile to improve the convergence rate by splitting methods~\cite{puterman2014markov}. Its fundamental architecture is to generate a sequence of value functions, $\{v^t\}$, $t=0,1,\ldots$ by an iterative process such that $v^t$ can sufficiently approach the robust team-optimal value function $v^*$ after a large number of iterations.
\par
The iterative process primarily encompasses four parts.
First, it starts with a given initial value function $v^0:=[v_{s^1}^0,v_{s^2}^0,\ldots,v_{s^m}^0]^T\in \mathcal{V}^0$ at $t=0$, where for ease of notations, we adopt $v_{s^j}^0$, $\forall j\in\{1,2,\ldots,m\}$ to represent the component of $v^0$ corresponding to state $s^j\in \mathcal{S}$ and $\mathcal{V}^0$ to represent the feasible set of $v^0$, respectively. Subsequently, the algorithm applies a process of policy improvement to obtain an improved team policy.
Specifically, given an arbitrary state $s^k\in\mathcal{S}$, for every joint action $a\in\mathcal{A}$, it first computes an approximation $\tilde{\rho}_{(s^k,a)}(v^t)$ of $\rho_{(s^k,a)}(v^t)$
with respect to $v^t:=[v_{s^1}^t,v_{s^2}^t,\ldots,v_{s^m}^t]^T$ at $t$, where $v_{s^j}^t$, $\forall j\in\{1,2,\ldots,m\}$ is the element of $v^t$ corresponding to state $s^j\in \mathcal{S}$, and $\rho_{(s^k,a)}(v^t)$ is the quantity obtained by running one step of the Gauss-Seidel value iteration~\cite{puterman2014markov} under the worst case of transition probability distributions, i.e.
\begin{equation}\label{eqn9}
\begin{split}
&\rho_{(s^k,a)}(v^t) =\min_{p(\cdot\mid s^k,a)\in \mathcal{P}_{a}(\cdot\mid k)}\left\{\sum_{l=1}^m r(s^k,a,s^l)p(s^l\mid s^k,a)\right. \\
&
 \left.+\lambda\left[\sum_{l<k}p(s^l\mid s^k,a)u^{t}_{0}(s^l)+\sum_{l\geq k}p(s^l\mid s^k,a)v^{t}_{s^l}\right]\right\},
\end{split}
\end{equation}
and $u_0^t(s^l) =
\max_{a\in\mathcal{A}}
\tilde{\rho}_{(s^l,a)}(v^t)$. One reason for using approximations therein is that an exact calculation is usually impractical when the state and/or action space is very large, and thus it requires to resort to approximation techniques. By selecting $d_{t+1}(s^k)\in \arg\max_{a\in\mathcal{A}}\tilde{\rho}_{(s^k,a)}(v^t)$ for every $s^k\in \mathcal{S}$,
one can then obtain an improved team policy $(d_{t+1})^\infty$, and subsequently each individual decision rule of $d_{t+1}$ is communicated to the corresponding player.
\par
Using the acquired intermediate value $u_0^t:=[u_0^t(s^1),u_0^t(s^2),\ldots,u_0^t(s^m)]^T$, the algorithm next implements a partial performance evaluation of the improved team policy $(d_{t+1})^\infty$ by running multiple steps of the value iteration similar to~(\ref{eqn9}). In this process, two new notations $\varsigma$ and $M_t$ of non-negative integers are introduced to represent the current iteration step and the total numbers of iterations, respectively. Specifically, for each given state $s^k\in\mathcal{S}$, every player first executes an action according to its received individual decision rule from $d_{t+1}$, and subsequently the algorithm collects the joint action $a'=d_{t+1}(s^k)\in\mathcal{A}$ of all players. Next, an approximation $\tilde{\rho}_{(s^k,a')}
(u^t_\varsigma)
$ of $\rho_{(s^k,a')}(u^t_\varsigma)$ is calculated with respect to $u^t_\varsigma:=[u_\varsigma^t(s^1),u_\varsigma^t(s^2),
\ldots,u_\varsigma^t(s^m)]^T$ at iteration $\varsigma$ $(\varsigma=0,1,\ldots,M_t)$. Therein, $\rho_{(s^k,a')}(u^t_\varsigma)$ is the quantity obtained by running one step of the Gauss-Seidel value iteration under the worst-case transition probability distribution used in~(\ref{eqn9}), i.e.
\begin{gather}\label{eqn6}
 \rho_{(s^k,a')}(u^t_\varsigma) =
\sum_{l=1}^m r(s^k,a',s^l)p^*(s^l\mid s^k,a')  \\
\nonumber +\lambda\left[\sum_{l<k}p^*(s^l\mid s^k,a')u^{t}_{\varsigma+1}(s^l)+\sum_{l\geq k}p^*(s^l\mid s^k,a')u^{t}_{\varsigma}(s^l)\right],
\end{gather}
and
$
u^t_{\varsigma+1}(s^l) =
\tilde{\rho}_{(s^l,a')}
(u^t_\varsigma)
$, where $
p^*(\cdot\mid s^k,a')
$ is the worst-case transition probability distribution used in~(\ref{eqn9}) for the given state $s^k$ and action $a'$. This iteration repeats by increment $\varsigma=\varsigma+1$ until $\varsigma=M_t$.
\par
Lastly, if the termination condition $\|u_0^t-v^t\|< (1-\lambda)\epsilon/2\lambda-\delta$ is satisfied at $t$, the algorithm stops and returns a team policy $(d^\epsilon)^\infty$, where $d^\epsilon=d_{t+1}$,
$\delta\in [0,(1-\lambda)^2\epsilon/2\lambda(1+\lambda))$ is the parameter used to confine the maximal approximation tolerance between $\rho_{(s,a)}(\cdot)$ and $\tilde{\rho}_{(s,a)}(\cdot)$ for $\forall (s,a)\in\mathcal{S}\times\mathcal{A}$, and $\epsilon>0$ is the small constant used to control the precision of the acquired $\epsilon$-robust team-optimal policy. The complete computational process is illustrated in Algorithm~1.
}

\begin{algorithm} \renewcommand{\algorithmicrequire}{\textbf{Input:}}
	\renewcommand{\algorithmicensure}{\textbf{Output:}}
	\caption{robust approximate
Team Policy Iteration (raTPI)}
	\label{alg:1}
	\begin{algorithmic}[1]
		\REQUIRE $v^0\in \mathcal{V}^0$, $\epsilon>0$, $\lambda\in [0,1)$, $\delta: 0\leq \delta<(1-\lambda)^2\epsilon/2\lambda(1+\lambda)$, and $\{M_t\}_{t\in \mathbb{N}}$.
        \STATE Set $t=0$.
		\STATE (Policy improvement) Set $k=1$ and go to 2(a).
        \begin{description}
          \item[(a)] For $\forall a\in \mathcal{A}$, first compute an approximation $\tilde{\rho}_{(s^k,a)}(v^t)$ of $\rho_{(s^k,a)}(v^t)$
               such that
              \begin{equation}\label{eqA1}
                |\tilde{\rho}_{(s^k,a)}(v^t)-\rho_{(s^k,a)}(v^t)|\leq \lambda \delta;
              \end{equation}
              then set
              \begin{equation}\label{eqA2}
                u_0^t(s^k)  =
                \max_{a\in\mathcal{A}}
              \tilde{\rho}_{(s^k,a)}(v^t),
              \end{equation}

              and select $d_{t+1}(s^k)\in \arg\max_{a\in\mathcal{A}}
              \tilde{\rho}_{(s^k,a)}(v^t)$.
          \item[(b)] If $k=m$, go to $3$; otherwise set $k=k+1$ and then return to 2(a).
        \end{description}
		\STATE (Partial policy evaluation)
        \begin{description}
          \item[(a)] If
          $
          \|u_0^t-v^t\|< (1-\lambda)\epsilon/2\lambda-\delta
          $ where $u_0^t:=[u_0^t(s^1),u_0^t(s^2),\ldots,u_0^t(s^m)]^T$,
              go to $4$; otherwise set $\varsigma=0$ and then go to 3(b).
          \item[(b)] If $\varsigma=M_t$, go to 3(f); otherwise set $k=1$ and then go to 3(c).
          \item[(c)] For $u_\varsigma^t:=[u_\varsigma^t(s^1),u_\varsigma^t(s^2),\ldots,u_\varsigma^t(s^m)]^T$, compute an approximation $\tilde{\rho}_{(s^k,d_{t+1}(s^k))}(u^t_\varsigma)$ of $\rho_{(s^k,d_{t+1}(s^k))}(u^t_\varsigma)$
              such that
              \begin{equation}\label{eqA3}
              |\tilde{\rho}_{(s^k,d_{t+1}(s^k))}(u^t_\varsigma)-
              \rho_{(s^k,d_{t+1}(s^k))}(u^t_\varsigma)|\leq \lambda \delta,
              \end{equation}
              and set
              \begin{equation}\label{eqA4}
              u^t_{\varsigma+1}(s^k) =
              \tilde{\rho}_{(s^k,d_{t+1}(s^k))}
              (u^t_\varsigma).
              \end{equation}

          \item[(d)] If $k=m$, go to 3(e); otherwise set $k=k+1$ and then return to 3(c).
          \item[(e)] Set $\varsigma=\varsigma+1$ and then return to 3(b).
          \item[(f)] Set $v^{t+1}=u_{M_t}^t$ and $t=t+1$, and then return to 2.
        \end{description}
		\STATE Set $d^\epsilon(s)=d_{t+1}(s)$ for $\forall s\in \mathcal{S}$ and then stop.
\ENSURE $\epsilon$-robust team-optimal policy $(d^{\epsilon})^\infty$.
	\end{algorithmic}
\end{algorithm}

\subsection{The convergence analysis}
{\color{red}Since directly proving the convergence of raTPI is difficult, we here divide the convergence analysis of raTPI into two parts.} The first part demonstrates the convergence of raTPI in a degenerated situation, i.e. $M_t=0$ for $\forall t\in \mathbb{N}$. {\color{red} Since in this case Algorithm~1 primarily relies on the value iteration to work and the partial policy evaluation will not come into effect, we refer to this degenerated form of raTPI as the robust approximate Team Value Iteration (raTVI).} Based on the results of raTVI, the second part establishes the convergence of raTPI in its general form where $\{M_t\}_{t\in \mathbb{N}}$ is an arbitrary sequence of non-negative integers.

\subsubsection{The degenerated form raTVI}
This part establishes the convergence of the degenerated form raTVI by setting $M_t=0$ for $\forall t\in \mathbb{N}$. We first consider the exact case where the tolerance of approximation computations is zero, i.e. $\delta=0$. {\color{red}Note that the partial policy evaluation in Algorithm~1 will not work due to $M_t=0$ for $\forall t\in \mathbb{N}$, and $\rho_{(s,a)}(\cdot) = \tilde{\rho}_{(s,a)}(\cdot)$ holds for $\forall (s,a)\in\mathcal{S}\times\mathcal{A}$ due to $\delta=0$. Thus, from~(\ref{eqA2}) and in view of $ v^{t+1}_s = u_{0}^t(s)$ for $\forall s\in \mathcal{S}$ in the step~3(f), one can get the iterative scheme of raTVI by}
\begin{equation}\label{eqn11}
v^{t+1}_s = \max_{a\in\mathcal{A}}\rho_{(s,a)}(v^t),\ \forall s\in\mathcal{S}.
\end{equation}
\par
{\color{red}
Utilizing the regular splitting in matrix iterative analysis~\cite{varga1962iterative}, we proceed to derive the vector form of~(\ref{eqn11}).
}
Let $\check{d}\in\mathcal{D}$ be the decision rule such that for $\forall s\in\mathcal{S}$, $\check{d}(s)\in \arg\max_{a\in\mathcal{A}}\rho_{(s,a)}(v^t)$, and $P_{\check{d}}\in \mathcal{P}_{\check{d}}$ be the worst-case transition probability matrix with element $P_{\check{d}}(l\mid k)$ in row $k$ and column $l$ when computing $\rho_{(s,\check{d}(s))}(v^t)$ by~(\ref{eqn9}). We then split $P_{\check{d}}$ into $P_{\check{d}}={P_{\check{d}}}^L+{P_{\check{d}}}^U$, {\color{red}
where ${P_{\check{d}}}^L$ and ${P_{\check{d}}}^U$ are the lower triangular matrix of $P_{\check{d}}$ and the upper triangular matrix including the diagonal elements of $P_{\check{d}}$, respectively.} As such, from~(\ref{eqn9}) and $ v^{t+1}_s = u_{0}^t(s)$ for $\forall s\in \mathcal{S}$ in the step~3(f), the vector form of~(\ref{eqn11}) can be given by
$
  v^{t+1}=r_{(\check{d},P_{\check{d}})}+\lambda\left[ {P_{\check{d}}}^L v^{t+1}+{P_{\check{d}}}^U v^{t}\right].
$
By re-arrangement, it yields
\begin{equation}\label{eqn10}
\begin{split}
    v^{t+1}&=(\underbrace{I-\lambda {P_{\check{d}}}^L}_{Q_{\check{d}}})^{-1} r_{(\check{d},P_{\check{d}})}+(I-\lambda {P_{\check{d}}}^L)^{-1}\underbrace{\lambda {P_{\check{d}}}^U}_{R_{\check{d}}}v^t \\
    &={Q_{\check{d}}}^{-1} r_{(\check{d},P_{\check{d}})}+{Q_{\check{d}}}^{-1}R_{\check{d}} v^t,
\end{split}
\end{equation}
where $(I-\lambda {P_{\check{d}}}^L)^{-1}$ exists because $\sigma(\lambda {P_{\check{d}}}^L)\leq \|\lambda {P_{\check{d}}}^L\|<1$. {\color{red}According to the definition of regular splitting (i.e. a pair of matrices $(B,C)$ is a regular splitting of matrix $A$ if $A=B-C$, $B^{-1}\geq 0$, and $C\geq 0$)~\cite{varga1962iterative}, one can see }that $(Q_{\check{d}},R_{\check{d}})$ is a regular splitting of $I-\lambda P_{\check{d}}$ because $Q_{\check{d}}-R_{\check{d}}=I-\lambda P_{\check{d}}$, ${Q_{\check{d}}}^{-1}=\sum_{\imath=0}^{\infty}(\lambda {P_{\check{d}}}^L)^\imath\geq 0$, and $R_{\check{d}}\geq 0$. For convenience, we refer to this specific regular splitting as the GS regular splitting afterwards, due to the use of the Gauss-Seidel (GS) value iteration. {\color{red}
That is, given a transition probability matrix $P$ and a constant $\lambda\in[0,1)$, the regular splitting $(Q,R)$ of $I-\lambda P$ is said to be a \emph{GS regular splitting} if $Q=I-\lambda P^L$ and $R=\lambda P^U$, where $P^L$ and $P^U$ are the lower triangular matrix of $P$ and the upper triangular matrix including the diagonal elements of $P$, respectively.
Since in~(\ref{eqn10}), $\check{d}\in\mathcal{D}$ is chosen such that the summation of the two terms on the right-hand side of~(\ref{eqn10}) is maximized with respect to the worse-case transition probability matrix $P_{\check{d}}\in\mathcal{P}_{\check{d}}$,
(\ref{eqn10}) is equivalent to
}
\begin{equation}\label{eqn13}
  v^{t+1}=\max_{d\in \mathcal{D}}\min_{P_d\in\mathcal{P}_d} \left\{Q_d^{-1}r_{(d,P_d)}+Q_d^{-1}R_d v^t\right\},
\end{equation}
where $(Q_d,R_d)$ is the GS regular splitting of $I-\lambda P_d$.
For any $v\in \mathcal{V}$, let the operator $\mathscr{Y}:\mathcal{V}\rightarrow \mathcal{V}$ be defined by
    \begin{equation}\label{eq14}
      \mathscr{Y} v:=\max_{d\in \mathcal{D}}\min_{P_d\in\mathcal{P}_d} \left\{Q_d^{-1}r_{(d,P_d)}+Q_d^{-1}R_d v\right\}.
    \end{equation}
Then, (\ref{eqn13}) can be rewritten by $v^{t+1}=\mathscr{Y} v^t$.
\par
Although the definition of $\mathscr{Y}$ is based on the GS regular splitting, we will show that whenever $(Q_d,R_d)$ is a regular splitting of $I-\lambda P_d$, $\mathscr{Y}$ is a contraction mapping on $\mathcal{V}$ and
the convergence rate of the sequence generated by it is less than $1$.
{\color{red}
Before formally presenting this result, we first introduce the definition of the convergence rate of a sequence~\cite{puterman2014markov}.
\begin{definition}\label{def3}
  Let $\{v^t\},\forall t\in \mathbb{N}$ be a sequence that converges to $v^*$. $\{v^t\}$ is said to converge at the order $c>0$ if there exists a constant $L$ such that $||v^{t+1}-v^*||\leq L||v^t-v^*||^c$ for $\forall t\in \mathbb{N}$, and its convergence rate is defined to be the smallest $L$ satisfying this inequality. The asymptotic average rate of convergence (AARC) of $\{v^t\}$ is defined as $\limsup_{t\rightarrow \infty}\left(
  ||v^{t}-v^*||/||v^{0}-v^*||\right)^{1/t}$. Given a non-negative function $f(t)\in \mathcal{V}$, $\{v^t\}$ is said to be $O(f(t))$ if $\limsup_{t\rightarrow \infty}
  ||v^{t}-v^*||/f(t)$ is finite. Any one of the above results is said to be global if it holds for $\forall v^0\in\mathcal{V}$; otherwise, it is said to be local.
\end{definition}

Having this definition, we now show in the following lemma that the sequence generated by $\mathscr{Y}$ will converge in norm to the robust team-optimal value function and its convergence rate is no greater than the constant $\alpha$ $(\alpha<1)$.
}
\begin{lemma}\label{lm1}
For any given $d\in \mathcal{D}$ and $P_d\in \mathcal{P}_d$, let $(Q_d,R_d)$ be a regular splitting of $I-\lambda P_d$. If $(Q_d,R_d)$ satisfies $\alpha:=\sup_{d\in \mathcal{D},P_d\in \mathcal{P}_d}\|Q_d^{-1}R_d\|<1$, then (a) the sequence $\{v^t\}$ generated by $v^{t+1}=\mathscr{Y} v^t$ converges in norm to the robust team-optimal value function $v^*$ for $t\rightarrow \infty$ and $v^*$ is the unique fixed point of $\mathscr{Y}$; and (b) $\{v^t\}$ converges globally at order $1$ at a rate no greater than $\alpha$; its global AARC is no greater than $\alpha$, and it is globally $O(\beta^t),\beta\leq \alpha$.

\end{lemma}
\par
{\color{red}
The detailed proof of this lemma is given in Appendix~\ref{app_A}.
%
From Lemma~1, one can see that only if the regular splitting $(Q_d,R_d)$ of $I-\lambda P_d$ satisfies $\alpha:=\sup_{d\in \mathcal{D},P_d\in \mathcal{P}_d}\|Q_d^{-1}R_d\|<1$, the sequence generated by the iterative scheme~(\ref{eqn13}) of raTVI will converge to the robust team-optimal value function. In reality, however, for the GS regular splitting, $\alpha<1$ will always hold. This fact is guaranteed by the following lemma.
\begin{lemma}\label{lm2}
(Proposition $6.3.5$ in~\cite{puterman2014markov})
   Given a transition probability matrix $P$, let $(Q_1,R_1)$ and $(Q_2,R_2)$ be the regular splittings of $I-\lambda P$ for $\lambda\in[0,1)$. If $R_2\leq R_1\leq\lambda P$, then $||Q_2^{-1}R_2||\leq||Q_1^{-1}R_1||$.
\end{lemma}
Note that for any given $d\in \mathcal{D}$ and $P_d\in \mathcal{P}_d$, the $(Q_d,R_d)$ in~(\ref{eqn13})
is the GS regular splitting of $I-\lambda P_d$, and $Q'_d=I$ and $R'_d=\lambda P_d$ are a trivial regular splitting of $I-\lambda P_d$. Moreover, since $R_d$ is the upper
triangular matrix including the diagonal elements of $\lambda P_d$ from the definition of the GS regular splitting, one can get $R_d\leq R'_d$. It then follows from Lemma~2 that $\|Q_d^{-1}R_d\|\leq \|{Q'_d}^{-1} R'_d\|=\lambda<1$. Since $d\in \mathcal{D}$ and $P_d\in \mathcal{P}_d$ are arbitrary, the condition
$\alpha<1$ for the $(Q_d,R_d)$ in (\ref{eqn13}) is always satisfied. From Lemma~1, one can then obtain the part~(a) of the following theorem.
}

\begin{theorem}\label{thm1}
For $M_t=0$, $\forall t\in \mathbb{N}$, if the approximation calculations are exact, i.e. $\delta=0$, in Algorithm~1, then for any $v^0\in \mathcal{V}$, (a) the sequence $\{v^t\}$ generated by the iterative scheme~(\ref{eqn13}) of raTVI converges in norm to the robust team-optimal value function $v^*$ for $t\rightarrow \infty$; it converges globally at order $1$ at a rate no greater than $\lambda$; its global AARC is no greater than $\lambda$, and it is globally $O(\beta^t),\beta\leq \lambda$; and (b) Algorithm~1 terminates within a finite number of iterations with an $\epsilon$-robust team-optimal policy $(d^{\epsilon})^\infty$ and its corresponding value function $v^{\epsilon}$ under the worst-case transition probability matrix satisfies $\|v^{\epsilon}-v^*\|< \epsilon$.
\end{theorem}

\begin{IEEEproof}
The result in part (a) is immediate from the above analysis.
We now prove part (b). Since $\{v^t\}$ converges to $v^*$, it is a Cauchy sequence. Hence, the termination condition in the step~3(a) of Algorithm~1 will be satisfied for any $\epsilon>0$ after a finite number of iterations in view of $v^{t+1}=u^t_0$ for $M_t=0$ in the step 3(f) of Algorithm~1. Without loss of generality, suppose that the algorithm terminates at $t=N$, i.e. $\|v^{N+1}-v^{N}\|< (1-\lambda)\epsilon/2\lambda$, and meanwhile it returns a decision rule $d^{\epsilon}$. Then, based on the setup in the step~2(a) and the step~4 of Algorithm~1, one can get from the iterative scheme~(\ref{eqn13}) that
$$
d^{\epsilon}\in \arg\max_{d\in \mathcal{D}} \left\{ \min_{P_d\in\mathcal{P}_d} \left[Q_d^{-1}r_{(d,P_d)}+Q_d^{-1}R_d v^{N}\right]\right\},
$$
where $(Q_d,R_d)$ is the GS regular splitting of $I-\lambda P_d$. Next, for the decision rule $d^{\epsilon}$, select
$$
P_{d^{\epsilon}}^{\epsilon}\in \arg\min_{P_{d^{\epsilon}}\in\mathcal{P}_{d^{\epsilon}}} \left[Q_{d^{\epsilon}}^{-1}r_{({d^{\epsilon}},P_{d^{\epsilon}})}
+Q_{d^{\epsilon}}^{-1}R_{d^{\epsilon}} v^{N}\right],
$$
where $(Q_{d^{\epsilon}},R_{d^{\epsilon}})$ is the GS regular splitting of $I-\lambda P_{d^{\epsilon}}$. Moreover, given $d\in \mathcal{D}$, $P_d\in\mathcal{P}_d$, and the GS regular splitting $(Q_d,R_d)$ of $I-\lambda P_{d}$ for $\lambda\in[0,1)$, we define the operator $\mathscr{T}_{(d,P_{d})}:\mathcal{V}\rightarrow \mathcal{V}$ for $v\in \mathcal{V}$ by
\begin{equation}\label{eqn15}
\mathscr{T}_{(d,P_{d})}v=
Q_{d}^{-1}r_{(d,
P_{d})}
+Q_{d}^{-1}R_{d} v.
\end{equation}
It is easy to check that $\mathscr{T}_{(d,P_{d})}$ is a contraction mapping on $\mathcal{V}$ because for any $u,v\in \mathcal{V}$,
\begin{equation*}
\begin{split}
&\|\mathscr{T}_{(d,P_{d})}u- \mathscr{T}_{(d,P_{d})}v\|
=\|Q_{d}^{-1}R_{d}
(u-v)\|
\\& \leq \|Q_{d}^{-1}R_{d}\|
 \|u-v\|
\leq \|I^{-1}(\lambda P_{d})\| \|u-v\|
=\lambda \|u-v\|,
\end{split}
\end{equation*}
where $\|Q_{d}^{-1}R_{d}\|
\leq \|I^{-1}(\lambda P_{d})\|$ follows from Lemma~\ref{lm2}. Then, by solving the fixed point equation $v=\mathscr{T}_{(d^{\epsilon},P_{d^{\epsilon}}^{\epsilon})}
v$ for variable $v\in\mathcal{V}$, one can get that the unique fixed point of $\mathscr{T}_{(d^{\epsilon},P_{d^{\epsilon}}^{\epsilon})}$
is $v_{(d^{\epsilon},P_{d^{\epsilon}}^{\epsilon})}$.  Also, from the definitions of $d^{\epsilon}$ and $P_{d^{\epsilon}}^{\epsilon}$, one can see that $\mathscr{T}_{(d^{\epsilon},P_{d^{\epsilon}}^{\epsilon})}v^{N}=\mathscr{Y} v^{N}=v^{N+1}$.
It follows that
\begin{equation*}
\begin{split}
&\|v_{(d^{\epsilon},P_{d^{\epsilon}}^{\epsilon})}-v^{N+1}\|
  =\|\mathscr{T}_{(d^{\epsilon},P_{d^{\epsilon}}^{\epsilon})}
  v_{(d^{\epsilon},P_{d^{\epsilon}}^{\epsilon})}
  -v^{N+1}\| \\
  &\leq \|\mathscr{T}_{(d^{\epsilon},P_{d^{\epsilon}}^{\epsilon})}
  v_{(d^{\epsilon},P_{d^{\epsilon}}^{\epsilon})}
  -\mathscr{T}_{(d^{\epsilon},P_{d^{\epsilon}}^{\epsilon})} v^{N+1}
  \|
  \\&+
  \|\mathscr{T}_{(d^{\epsilon},P_{d^{\epsilon}}^{\epsilon})} v^{N+1}-\mathscr{T}_{(d^{\epsilon},P_{d^{\epsilon}}^{\epsilon})}v^{N} \|
  \\&\leq \lambda\|v_{(d^{\epsilon},P_{d^{\epsilon}}^{\epsilon})}-v^{N+1}\|
  +\lambda \|v^{N+1}-v^N\|.
\end{split}
\end{equation*}
By re-arrangement, one can get $\|v_{(d^{\epsilon},P_{d^{\epsilon}}^{\epsilon})}-v^{N+1}\|
\leq \lambda\|v^{N+1}-v^{N}\|/(1-\lambda)<\epsilon/2$.
Moreover, by applying the contraction property of $\mathscr{Y}$, we have
\begin{equation*}
\begin{split}
&\|v^{N+1}- v^{*}\| \leq \|v^{N+1}-\mathscr{Y}v^{N+1}\| + \|\mathscr{Y}v^{N+1}- v^{*}\| \\&
\leq \lambda\|v^{N+1} - v^{N}\|+ \lambda \|v^{N+1}- v^{*}\|.
\end{split}
\end{equation*}
By re-arrangement,
 $\|v^{N+1}-v^*\|
\leq \lambda\|v^{N+1}-v^{N}\|/(1-\lambda)< \epsilon/2$. Consequently,
$
 \|v_{(d^{\epsilon},P_{d^{\epsilon}}^{\epsilon})}-v^{*}\|
 \leq \|v_{(d^{\epsilon},P_{d^{\epsilon}}^{\epsilon})}-v^{N+1}\|
+
\|v^{N+1}-v^*\|  < \epsilon.
$
Since $v_{(d^{\epsilon},P_{d^{\epsilon}}^{\epsilon})}$
is the value function corresponding to the policy
$(d^\epsilon)^\infty$ under the worst-case transition probability matrix, we have the result in part~(b).
\end{IEEEproof}

{\color{red}
One important precondition for the establishment of Theorem~1 is that those approximation calculations for $\rho_{(s,a)}(\cdot)$, $\forall (s,a)\in \mathcal{S}\times\mathcal{A}$ in Algorithm~1 are exact, i.e. $\delta=0$. In contrast, if the error is inevitable for approximations, the establishment of the convergence will become quite intractable. This is because the approximation error will be continuously accumulated along with the iteration of the algorithm and eventually it will lead the iterative sequence to deviate from the desired value function. Fortunately, however, we find that if the parameter $\delta$ is bounded by $ 0\leq \delta<(1-\lambda)^2\epsilon/2\lambda(1+\lambda)$, every element $ \tilde{v}^t$ of the sequence generated by Algorithm~1 will
always lie within a small interval, which is centered at the element $v^t$ of the sequence generated by Algorithm~1 for $\delta=0$. That is, for the sequence $ \{\tilde{v}^t\}$ generated by raTVI with $  \delta\in [0,(1-\lambda)^2\epsilon/2\lambda(1+\lambda))$ and the sequence $\{v^t\}$ generated by raTVI with $\delta=0$, the following lemma holds.
}

\begin{lemma}\label{lm3}
  If $v^0=\tilde{v}^0$, then $\{v^t\}$ and $\{\tilde{v}^t\}$ satisfy
  \begin{equation}\label{eqn14}
  v^t-\theta\mathbf{1}\leq \tilde{v}^t \leq v^t+\theta\mathbf{1},\ \forall t\in \mathbb{N},
  \end{equation}
  where $\theta:={\lambda \delta}/(1-\lambda)\in [0,\epsilon/2)$.
\end{lemma}
{\color{red}
See Appendix~\ref{app_B} for its proof. From this lemma, one can immediately get $\lim_{t\rightarrow \infty}\tilde{v}^t=v^*$ by noticing that $\lim_{t\rightarrow \infty}v^t=v^*$ from Theorem~1 and $\theta\in [0,\epsilon/2)$. We display this result in the following theorem.
}
\begin{theorem}\label{thm2}
For $M_t=0$, $\forall t\in \mathbb{N}$, if the tolerance of approximation calculations is within $\lambda\delta$ for $\delta\in [0,(1-\lambda)^2\epsilon/2\lambda(1+\lambda))$ in Algorithm~1, then for any $v^0\in \mathcal{V}$, (a) the sequence $\{\tilde{v}^t\}$ generated by Algorithm~1 with the initialization $\tilde{v}^0=v^0$ converges in norm to the robust team-optimal value function $v^*$ for $t\rightarrow \infty$; it converges globally at order $1$ at a rate no greater than $\lambda$; its global AARC is no greater than $\lambda$, and it is globally $O(\beta^t),\beta\leq \lambda$; and (b) those results shown in the part~(b) of Theorem~\ref{thm1} still hold.
\end{theorem}

\begin{IEEEproof}
Since we have got in Theorem~\ref{thm1} that the sequence $\{v^t\}$ generated by Algorithm~1 for $M_t=0$, $\forall t\in \mathbb{N}$ and $\delta=0$ converges to $v^*$, there is a positive integer $N$ for any $\epsilon>0$ such that $\|v^t-v^*\|< \epsilon /2$ holds for $\forall t\geq N$. Moreover, from Lemma~\ref{lm3}, one can obtain that for $\forall t\in \mathbb{N}$,
\begin{equation}\label{eqn19}
\|\tilde{v}^t-v^*\|\leq \|\tilde{v}^t-v^t\| + \|v^t-v^*\|\leq \|v^t-v^*\|+\theta.
\end{equation}
Therefore, we have $\|\tilde{v}^t-v^*\|\leq \epsilon/2 + \theta<\epsilon$ for $\forall t\geq N$, which implies that $\tilde{v}^t$ converges to $v^*$.
\par
We now calculate the convergence rates of $\{\tilde{v}^t\}$. From Theorem~\ref{thm1}, we first note that $\{v^t\}$ satisfies $\|v^t-v^*\|\leq \lambda \|v^{t-1}-v^*\|$ for $\forall t\in \mathbb{N}$. Moreover, from Lemma~\ref{lm3}, $\|v^{t-1}-v^*\|\leq \| v^{t-1}- \tilde{v}^{t-1}\| + \|\tilde{v}^{t-1}-v^*\|\leq \|\tilde{v}^{t-1}-v^*\|+\theta$ holds for all $t$. Based on these two inequalities, it follows from~(\ref{eqn19}) that
\begin{equation*}
\begin{split}
&\|\tilde{v}^t-v^*\| \leq \|v^t-v^*\|+\theta
\\
&
  \leq \lambda \|v^{t-1}-v^*\|+\theta
  \leq \lambda \|\tilde{v}^{t-1}-v^*\|+(\lambda+1)\theta
.
\end{split}
\end{equation*}
Since $\theta\in[0,\epsilon/2)$ can be arbitrarily close to zero in view of the arbitrariness of $\epsilon$, $\|\tilde{v}^t-v^*\| \leq \lambda \|\tilde{v}^{t-1}-v^*\|$ for all $t$. It means that $\{\tilde{v}^t\}$ converges globally at order $1$ and the convergence rate is no greater than $\lambda$. Similarly, from~(\ref{eqn19}), we also have $\|\tilde{v}^t-v^*\|\leq \|v^t-v^*\|$. It follows that
$
  {\|\tilde{v}^t-v^*\|}/{\|v^0-v^*\|}\leq
  {\|v^t-v^*\|}/{\|v^0-v^*\|}$
  and
   ${\|\tilde{v}^t-v^*\|}/{\beta^t} \leq {\|v^t-v^*\|}/{\beta^t}
  $  for $\beta>0$ and $\forall t\in \mathbb{N}$.
Therefore,
based on the results in Theorem~\ref{thm1}, one can get that the global AARC of $\{\tilde{v}^t\}$ is no greater than $\lambda$ and it is globally $O(\beta^t),\beta\leq \lambda$.
\par
We next prove part~(b). Since $\{\tilde{v}^t\}$ converges to $v^*$ from the result in part~(a), the termination condition in the step~3(a) of Algorithm~1 will be satisfied after a finite number of iterations.
Moreover, since when $\delta=0$, it is shown from Theorem~\ref{thm1}(b) that $\{v^t\}$ will satisfy the termination condition within a finite number of iterations, there exists a positive integer $N$ for any $\epsilon_1>0$ such that $\|v^{N+1}-v^N\|<(1-\lambda)\epsilon_1/2\lambda$. Then, from Lemma~\ref{lm3}, we have
\begin{equation*}
\begin{split}
&\|\tilde{v}^{N+1}-\tilde{v}^{N}\| \leq  \|\tilde{v}^{N+1}-  v^{N+1} \|+
\|v^{N+1}-v^{N}\|
\\
&
+\|v^{N}-\tilde{v}^{N}\|
\leq
  \|v^{N+1}-v^{N}\|+2\theta
\\
&
< (1-\lambda)\epsilon_1/2\lambda
+2\lambda \delta/(1-\lambda)
.
\end{split}
\end{equation*}
Given that $\epsilon_1$ is arbitrary, we therefore select a specific $\epsilon_1$ such that $\epsilon_1\leq \epsilon -2\lambda(1+\lambda)\delta/(1-\lambda)^2$ for any given $\epsilon>0$. As such, $\|\tilde{v}^{N+1}-\tilde{v}^{N}\| \leq (1-\lambda)\epsilon/2\lambda-\delta$ holds. It implies that for $M_t=0$, $\forall t\in \mathbb{N}$, Algorithm~1 with approximation errors will also terminate at $t=N$. Moreover, for any given $s\in \mathcal{S}$, one can find that there exists a constant gap between $\rho_{(s,a)}(v^t)$ and $\rho_{(s,a)}(\tilde{v}^t)$ (see the proof of Lemma~\ref{lm3}) and an approximation error between $\rho_{(s,a)}(\tilde{v}^t)$ and $\tilde{\rho}_{(s,a)}(\tilde{v}^t)$ from~(\ref{eqA1}) for $\forall a \in \mathcal{A}$ and $\forall t\in \mathbb{N}$, and they do not affect the selection of the maximizing actions. Thus,
$\arg\max_{a\in\mathcal{A}}\rho_{(s,a)}(v^N)=\arg\max_{a\in\mathcal{A}}\tilde{\rho}_{(s,a)}(\tilde{v}^N)$.
Without loss of generality, we choose the same decision rule for Algorithm~1 in the case with and without approximation errors at $t=N$. Thus, the result $\|v^{\epsilon}-v^*\|< \epsilon$ shown in Theorem~1(b) also holds for $\delta\in [0,(1-\lambda)^2\epsilon/2\lambda(1+\lambda))$.

%
%

\end{IEEEproof}

\subsubsection{The general form}
{\color{red}
We have presented the convergence of the degenerated form raTVI by the aforementioned analysis. Based on the results in this specific case, we now proceed to demonstrate the convergence of raTPI in its general form, i.e. $M_t$ is a non-negative integer for $\forall t\in \mathbb{N}$.
}
Via a similar analytical process, we first give the results in the exact case where the approximation calculations in Algorithm~1 can be exactly obtained, i.e. $\delta=0$, and then extend them to the inexact case where the tolerance of approximation calculations is within $\lambda \delta$ for $ 0\leq \delta<(1-\lambda)^2\epsilon/2\lambda(1+\lambda)$.
\par

We begin with considering the exact case, i.e. $\delta=0$. Let the operator $\mathscr{B}: \mathcal{V}\rightarrow \mathcal{V}$ for $v\in \mathcal{V}$ be defined by
\begin{equation}\label{eq18}
\begin{split}
\mathscr{B}v
&:=\mathscr{Y} v - v
\\&
=\max_{d\in \mathcal{D}}\min_{P_d\in \mathcal{P}_d}
\left\{ Q_d^{-1}r_{(d,P_d)} + Q_d^{-1}R_d v - v\right\},
\end{split}
\end{equation}
where $(Q_d,R_d)$ is the GS regular splitting of $I-\lambda P_d$. From this definition, it is easy to see that the fixed-point of $\mathscr{Y}$ is the same as the zero-point of $\mathscr{B}$. Moreover, for any given $v\in \mathcal{V}$, we denote the set of \emph{$v$-improving decision rules} by
$$
\mathcal{D}_v:=\arg\max_{d\in \mathcal{D}}\left\{\min_{P_d\in \mathcal{P}_d}
\left[
Q_d^{-1}r_{(d,P_d)}+
Q_d^{-1}R_d v \right]\right\},
$$
and accordingly denote the set of \emph{$d_v$-decreasing transition probability matrices} for a given $d_v\in \mathcal{D}_v$ by
$$
\mathcal{P}_{d_v}^* := \arg\min_{P_{d_v}\in \mathcal{P}_{d_v}}
\left[
{Q_{d_v}}^{-1}r_{(d_v,P_{d_v})}+
{Q_{d_v}}^{-1}R_{d_v} v \right].
$$
{\color{red}
Having these concepts at hand, one can immediately obtain the following proposition.}
\begin{proposition}\label{prop1}
  For any given $u,v\in \mathcal{V}$ and $d_v\in \mathcal{D}_v$, there exists a $P'_{d_v}\in \mathcal{P}_{d_v}$ such that
  $
  \mathscr{B}u\geq \mathscr{B}v +({Q'_{d_v}}^{-1}R'_{d_v}-I)(u-v)
  $,
   where $(Q'_{d_v},R'_{d_v})$ is the GS regular splitting of $I-\lambda P'_{d_v}$.
\end{proposition}
{\color{red}
See Appendix~\ref{app_C} for its proof. Intuitively, this proposition implicitly suggests that $\mathscr{B}$ is a Lipschitz continuous mapping by selecting appropriate norms.
}
\par
On the other hand, one can conveniently derive the vector formulation of the iterative scheme of raTPI. For any $v^t$, $\forall t\in \mathbb{N}$,
let $d_{v^t}$ be the decision rule satisfying $d_{v^t}(s) \in \arg\max_{a\in\mathcal{A}}\rho_{(s,a)}(v^t)$ for $\forall s\in \mathcal{S}$, and $P_{d_{v^t}}^*$ be the transition probability matrix with element
$
p^*(s^l \mid s^k,d_{v^t}(s^k))
$ in row $k$ and column $l$ for $ s^k, s^l \in \mathcal{S}$, where $
p^*(\cdot \mid s^k,d_{v^t}(s^k))
$ is the worst-case transition probability distribution used to calculate $\rho_{(s^k,d_{v^t}(s^k))}(v^t)$ by~(\ref{eqn9}).
Then, from (\ref{eqA1}) and (\ref{eqA2}), by applying an analogous argument to the derivation of~(\ref{eqn10}), one can get the vector form of (\ref{eqA2}) for $\delta=0$ by
\begin{equation}\label{eqn21}
\begin{split}
    u_0^t&=(\underbrace{I-\lambda {P_{d_{v^t}}^*}^L}_{Q_{d_{v^t}}^*})^{-1}
    r_{(d_{v^t},P_{d_{v^t}}^*)}+(I-\lambda {P_{d_{v^t}}^*}^L)^{-1}\underbrace{\lambda {P_{d_{v^t}}^*}^U}_{R_{d_{v^t}}^*}v^t \\
    &={Q_{d_{v^t}}^*}^{-1}r_{(d_{v^t},P_{d_{v^t}}^*)}
    +{Q_{d_{v^t}}^*}^{-1} R_{d_{v^t}}^* v^t,
\end{split}
\end{equation}
where
${P_{d_{v^t}}^*}^L$ and ${P_{d_{v^t}}^*}^U$
are the lower triangular matrix of $P_{d_{v^t}}^*$ and the upper triangular matrix including the diagonal elements of $P_{d_{v^t}}^*$, respectively.
Clearly, one can see that $(Q_{d_{v^t}}^*,R_{d_{v^t}}^*)$ is the GS regular splitting of $I-\lambda P_{d_{v^t}}^*$, $d_{v^t}\in \mathcal{D}_{v^t}$, and $P_{d_{v^t}}^*\in \mathcal{P}_{d_{v^t}}^*$ from their definitions.
{\color{red}
Moreover, based on the definition of $\mathscr{T}_{(d,P_{d})}$ in~(\ref{eqn15}), one can rewrite (\ref{eqn21}) by $u_0^t=\mathscr{T}_{(d_{v^t},P_{d_{v^t}}^*)}v^t$. Also, via an analogous derivation from~(\ref{eqn6}), (\ref{eqA3}), and (\ref{eqA4}) for $\delta=0$, one can get the vector form of (\ref{eqA4}) by $u^t_{\varsigma +1}= \mathscr{T}_{(d_{v^t},P_{d_{v^t}}^*)}u^t_{\varsigma }$. In view of $v^{t+1}=u_{M_t}^t$ in the step~3(f) of Algorithm~1,
one can then obtain the iterative scheme of raTPI for $\delta=0$ by
}
\begin{equation}\label{eqn22}
\begin{split}
  v^{t+1}&=\left(
  \mathscr{T}_{(d_{v^t},P_{d_{v^t}}^*)}
  \right)^{M_t+1}
  v^t \\
  &=\sum_{\varsigma=0}^{M_t}\left({Q_{d_{v^t}}^*}^{-1} R_{d_{v^t}}^*\right)^{\varsigma}{Q_{d_{v^t}}^*}^{-1}
  r_{(d_{v^t},P_{d_{v^t}}^*)}
  \\
  &+\left({Q_{d_{v^t}}^*}^{-1} R_{d_{v^t}}^*\right)^{M_t+1}v^t  \\
  &= v^t + \sum_{\varsigma=0}^{M_t}\left({Q_{d_{v^t}}^*}^{-1} R_{d_{v^t}}^*\right)^{\varsigma} \left(
  {Q_{d_{v^t}}^*}^{-1} r_{(d_{v^t},P_{d_{v^t}}^*)} \right.
  \\
  &
  \left. +
  {Q_{d_{v^t}}^*}^{-1} R_{d_{v^t}}^* v^t -v^t \right) \\
  &= v^t + \sum_{\varsigma=0}^{M_t}\left({Q_{d_{v^t}}^*}^{-1} R_{d_{v^t}}^*\right)^{\varsigma}\left( \mathscr{B}v^t \right),
\end{split}
\end{equation}
{\color{red}
where the second equality follows from the definition of $\mathscr{T}_{(d_{v^t},P_{d_{v^t}}^*)}$
in~(\ref{eqn15}) and the fourth equality follows from~(\ref{eq18}) in view of $d_{v^t}\in \mathcal{D}_{v^t}$ and $P_{d_{v^t}}^*\in \mathcal{P}_{d_{v^t}}^*$.
}
\par
From~(\ref{eqn22}), one can see that raTPI incorporates raTVI as a special case, because when $M_t=0$ for $\forall t\in \mathbb{N}$, (\ref{eqn22}) will degenerate to the iterative scheme of raTVI, $v^{t+1}=\mathscr{Y} v^t$. Moreover, for another extreme case $M_t\rightarrow \infty$, (\ref{eqn22}) will reduce to the performance evaluation of the improved decision rule $d_{v^t}\in \mathcal{D}_{v^t}$ for $v^t$ under the worst-case transition probability matrix $P_{d_{v^t}}^*\in \mathcal{P}_{d_{v^t}}^*$, i.e.
\begin{equation}\label{eqn23}
\begin{split}
  v^{t+1}&= v^t + \sum_{\varsigma=0}^{\infty}\left({Q_{d_{v^t}}^*}^{-1} R_{d_{v^t}}^*\right)^{\varsigma}( \mathscr{B}v^t )  \\
  &=v^t + \left(I-{Q_{d_{v^t}}^*}^{-1} R_{d_{v^t}}^*\right)^{-1}
  \left[{Q_{d_{v^t}}^*}^{-1}r_{(d_{v^t},P_{d_{v^t}}^*)} \right.
  \\
  &
  \left.-\left(I-{Q_{d_{v^t}}^*}^{-1} R_{d_{v^t}}^*\right)v^t
  \right]  \\
  &= \left(I-{Q_{d_{v^t}}^*}^{-1} R_{d_{v^t}}^*\right)^{-1}
  {Q_{d_{v^t}}^*}^{-1}  r_{(d_{v^t},P_{d_{v^t}}^*)} \\
  &= \left(I-\lambda P_{d_{v^t}}^*\right)^{-1} r_{(d_{v^t},P_{d_{v^t}}^*)}
  =v_{(d_{v^t},P_{d_{v^t}}^*)},
\end{split}
\end{equation}
{\color{red}
where the second equality follows from the third equality in~(\ref{eqn22}) by letting $M_t\rightarrow \infty$ and using $\sum_{\varsigma=0}^{\infty}({Q_{d_{v^t}}^*}^{-1} R_{d_{v^t}}^*)^{\varsigma}
=
(I-{Q_{d_{v^t}}^*}^{-1} R_{d_{v^t}}^*)^{-1}
$.
}
\par
{\color{red}
To analyze the convergence of the
iterative scheme~(\ref{eqn22}), we construct two auxiliary operators.
}
Given a non-negative integer $M\in \mathbb{N}$, for $v\in \mathcal{V}$, we define the operator $\mathscr{W}^M: \mathcal{V}\rightarrow \mathcal{V}$ by
\begin{equation}\label{eq22n}
  \mathscr{W}^{M} v := (
  \mathscr{T}_{(d_{v},P_{d_{v}}^*)}
  )^{M+1}
  v,
\end{equation}
and the operator $\mathscr{U}^M: \mathcal{V}\rightarrow \mathcal{V}$ by
\begin{equation}\label{eq23}
  \mathscr{U}^M v := \max_{d\in \mathcal{D}}\max_{P_d\in \mathcal{P}_d} \Phi(d,P_d,v),
\end{equation}
where $d_v\in\mathcal{D}_v$, $P_{d_{v}}^*\in \mathcal{P}_{d_{v}}^*$, and
$
\Phi(d,P_d,v)
:=
  \sum_{\varsigma=0}^{M}\left(Q_{d}^{-1} R_{d}\right)^{\varsigma}Q_{d}^{-1}
  r_{(d,P_{d})}
  +\left(Q_{d}^{-1} R_{d}\right)^{M+1}v
$ for the GS regular splitting $(Q_d,R_d)$ of $I-\lambda P_d$. In particular, based on~(\ref{eq22n}), the iterative scheme (\ref{eqn22}) of raTPI for $\delta=0$ can be rewritten by $v^{t+1}=\mathscr{W}^{M_t} v^t$. In the following, we show some properties of these two operators, which will become the basis of proving the convergence of raTPI.

\begin{lemma}\label{lm4}
For the operator $\mathscr{U}^M$, (a) it is a contraction mapping on $\mathcal{V}$ with constant $\lambda^{M+1}$; and (b) the sequence $\{\omega^t\}$ generated by $\omega^{t+1}=\mathscr{U}^M\omega^t$ for any $\omega^0\in \mathcal{V}$ converges in norm to the robust team-optimal value function $v^*$, which is the unique fixed point of $\mathscr{U}^M$.
\end{lemma}


{\color{red}
This lemma (see Appendix~\ref{app_D} for its proof) guarantees that the sequence generated by $\mathscr{U}^M$ will converge to the robust team-optimal value function $v^*$ with a convergence rate at least $\lambda^{M+1}$. Also, note from Lemma~\ref{lm1} that
$v^*$ is the limit point of convergence for the sequence generated by $\mathscr{Y}$. Therefore, it motivates us to postulate the following fact that for any $M\in\mathbb{N}$ and $\forall t\in \mathbb{N}$,
\begin{equation}\label{eq22}
(\mathscr{U}^M)^t v^0\geq (\mathscr{W}^M)^t v^0\geq
(\mathscr{Y})^t v^0
\end{equation}
holds for some initial value functions $v^0$. This is because once such a presumption holds, the sequence generated by $\mathscr{W}^M$ (i.e. the iterative scheme~(\ref{eqn22}) of raTPI) will also converge to $v^*$ for $t\rightarrow\infty$. As such, the convergence of raTPI for $\delta=0$ to the robust team-optimal value function will be established. We show by the following two lemmas that such a postulation is indeed sound for $v^0\in
\mathcal{V}_{\mathscr{B}}:=\{v\in \mathcal{V} \mid \mathscr{B}v\geq 0\}
$.
}

\begin{lemma}\label{lm5}
  For any $u,v\in \mathcal{V}$ satisfying $u\geq v$, $\mathscr{U}^M u \geq \mathscr{W}^M v$ holds for any $M\in \mathbb{N}$. Moreover, if $u\in \mathcal{V}_{\mathscr{B}}$, then $\mathscr{W}^M u\geq \mathscr{Y} v$.
\end{lemma}

\begin{lemma}\label{lm6}
  If $v\in \mathcal{V}_{\mathscr{B}}$, then $\mathscr{W}^M v\in \mathcal{V}_{\mathscr{B}}$ for any $M \in \mathbb{N}$.
\end{lemma}
{\color{red}
See Appendix~\ref{app_E} and~\ref{app_F} for their proofs, respectively. Based on the above three lemmas, one can then get the convergence result of raTPI for $\delta=0$, which is shown in the following theorem.
}

\begin{theorem}\label{thm3}
For any non-negative integer sequence $\{M_t\}_{t\in \mathbb{N}}$, if the approximation calculations are exact, i.e. $\delta=0$, in Algorithm~1, then for any $v^0\in \mathcal{V}_{\mathscr{B}}$,
(a)~the sequence $\{v^t\}$ generated by the iterative scheme (\ref{eqn22}) of raTPI converges monotonically and in norm to the robust team-optimal value function $v^*$;
(b)~those results shown in the part~(b) of Theorem~\ref{thm1} still hold; and
(c)~let $d_{v^t}$ and $d_{v^*}$ be the $v^t$-improving and $v^*$-improving decision rules, respectively. Then, there exists a $P_{d_{v^t}}^{*}\in \mathcal{P}_{d_{v^t}}^*$ and a $P'_{d_{v^*}}\in \mathcal{P}_{d_{v^*}}$ for which
\begin{equation*}
\begin{split}
        &\|v^{t+1}-v^*\|\leq \| v^t - v^*  \|
        \\&
        \times
        \left(
     \| {Q_{d_{v^t}}^*}^{-1} R_{d_{v^t}}^* -
     {Q'_{d_{v^*}}}^{-1} R'_{d_{v^*}} \|\frac{1-\lambda^{M_t}}{1-\lambda} + \lambda^{M_t+1}
     \right) ,
\end{split}
\end{equation*}
        where $(Q_{d_{v^t}}^*,R_{d_{v^t}}^*)$ and $(Q'_{d_{v^*}},R'_{d_{v^*}})$ are the GS regular splitting of $I-\lambda P_{d_{v^t}}^{*}$ and $I-\lambda P'_{d_{v^*}}$, respectively.
        In particular, if
  $\lim_{t\rightarrow \infty}
        \| P_{d_{v^t}}^{*} - P'_{d_{v^*}} \|=0$, then there exists a $K\in \mathbb{N}$ for any $\epsilon >0$ such that
     $
     \|v^{t+1} - v^* \| \leq
     \left(\epsilon + \lambda^{M_t+1}\right)\| v^t - v^*  \|
     $
     for $t\geq K$.
\end{theorem}

\begin{IEEEproof}
We first prove part~(a). Let  $\{y^t\}$ and $\{\omega^t\}$ be the sequences generated by $y^{t+1}=\mathscr{Y} y^t$ and $\omega^{t+1}=\mathscr{U}^{M_t}\omega^t$, respectively, where $y^0=\omega^0=v^0$. In the following, we show by induction that $v^t\in \mathcal{V}_{\mathscr{B}}$, $v^{t+1}\geq v^{t}$, and $\omega^t\geq v^t \geq y^t$ for $\forall t\in \mathbb{N}$.
\par
First, when $t=0$, $v^0 \in \mathcal{V}_{\mathscr{B}}$ and $y^0=\omega^0=v^0$ follow from the assumption. Moreover, from Lemma~\ref{lm5}, we have $v^1=\mathscr{W}^{M_0} v^0 \geq \mathscr{Y} v^0\geq v^0$. Consequently, the induction hypothesis holds for $t=0$. Suppose now that it is satisfied for $t=\kappa$. Then, applying Lemma~\ref{lm6} leads to $v^{\kappa+1}=\mathscr{W}^{M_\kappa}v^\kappa\in \mathcal{V}_{\mathscr{B}}$. From the definition of $\mathscr{W}^{M}$ in~(\ref{eq22n}) and the fourth equality in~(\ref{eqn22}),
moreover, $v^{\kappa+2}=\mathscr{W}^{M_{\kappa+1}}v^{\kappa+1}$ can further be given by
$$
  v^{\kappa+2}=v^{\kappa+1} + \sum_{\varsigma=0}^{M_{\kappa+1}}\left({Q_{d_{v^{\kappa+1}}}^*}^{-1} R_{d_{v^{\kappa+1}}}^*\right)^{\varsigma}\left( \mathscr{B}v^{\kappa+1} \right) \geq v^{\kappa+1}.
$$
Since $\omega^\kappa \geq v^\kappa \geq y^\kappa$ from the hypothesis, it follows from Lemma~\ref{lm5} that
$
\omega^{\kappa+1}=\mathscr{U}^{M_\kappa} \omega^{\kappa}\geq \mathscr{W}^{M_\kappa} v^\kappa=v^{\kappa+1} \geq \mathscr{Y} y^{\kappa}=y^{\kappa+1}
$.
Thus, the induction hypothesis is satisfied for $t=\kappa+1$. Note from Lemma~\ref{lm1} and Lemma~\ref{lm4} that both $y^t$ and $\omega^t$ converge in norm to $v^*$ for $t\rightarrow \infty$. Therefore, we have $v^*\leq\lim_{t\rightarrow\infty}v^t\leq v^*$, which means that $v^t$ will converge in norm to $v^*$. Part (a) is established.
\par
We proceed to prove part (b). Since part (a) has shown that the sequence $\{v^t\}$ generated by the iterative scheme of raTPI for $\delta=0$
is convergent, the termination condition in the step~3(a) of Algorithm~1 will be satisfied after a finite number of iterations by noticing
\begin{equation*}
  \begin{split}
  v^{t+1}&= (\mathscr{T}_{(d_{v^t},P_{d_{v^t}}^*)}
  )^{M_t+1}
  v^t
  =\mathscr{W}^{M_t} v^t
  \\ &
  \geq u^t_0=\mathscr{T}_{(d_{v^t},P_{d_{v^t}}^*)}v^t
  =\mathscr{Y}v^t\geq v^t,
  \end{split}
  \end{equation*}
  {\color{red}
  where the first line follows from~(\ref{eqn22}) and~(\ref{eq22n}); the first inequality follows from Lemma~\ref{lm5} and (\ref{eqn21}); and the second inequality follows from $v^t\in \mathcal{V}_{\mathscr{B}}$.
  }
  Suppose now that Algorithm~1 terminates at $t=N$ and returns a policy $(d^\epsilon)^\infty$. From the step~2(a) and~4 in Algorithm~1, it is known that $d^\epsilon$ is the $v^N$-improving decision rule. Let $P_{d^\epsilon}^*\in \mathcal{P}_{d^\epsilon}^*$ be the $d^\epsilon$-decreasing transition probability matrix with respect to $v^N$. Then, from the definition of $\mathscr{Y}$, we have $\mathscr{T}_{(d^\epsilon,P_{d^\epsilon}^*)}v^N
  =\mathscr{Y}v^N$. Moreover, by applying the contraction property of $\mathscr{Y}$, one can get
  \begin{equation*}
  \begin{split}
  &\|\mathscr{T}_{(d^\epsilon,P_{d^\epsilon}^*)}v^N-v^*\|
  =\|\mathscr{Y}v^N-v^*\|
  \leq \|\mathscr{Y}v^N -(\mathscr{Y})^2 v^N\|
  \\&
  +\|(\mathscr{Y})^2 v^N - v^*\|
  \leq \lambda \|v^N -\mathscr{Y} v^N\| + \lambda \|\mathscr{Y} v^N - v^*\|
  \\&
  =\lambda \|v^N -\mathscr{T}_{(d^\epsilon,P_{d^\epsilon}^*)}v^N\| + \lambda \|\mathscr{T}_{(d^\epsilon,P_{d^\epsilon}^*)}v^N - v^*\|.
  \end{split}
  \end{equation*}
  Since Algorithm~1 terminates at $t=N$, i.e. $\|u_0^N-v^N\|=
  \|\mathscr{T}_{(d^\epsilon,P_{d^\epsilon}^*)}v^N-v^N\|
  <(1-\lambda)\epsilon/2\lambda$, by re-arrangement the above inequality can lead to
  \begin{equation}\label{eq25}
  \|\mathscr{T}_{(d^\epsilon,P_{d^\epsilon}^*)}v^N-v^*\|
  \leq \|v^N -\mathscr{T}_{(d^\epsilon,P_{d^\epsilon}^*)}v^N\|
  \frac{\lambda}
  {1-\lambda}
  <\epsilon/2.
  \end{equation}
  On the other hand, since
  $
  (\mathscr{T}_{(d^\epsilon,P_{d^\epsilon}^*)})^\infty v^N- v^N=(\mathscr{T}_{(d^\epsilon,P_{d^\epsilon}^*)} v^N-v^N)
  +((\mathscr{T}_{(d^\epsilon,P_{d^\epsilon}^*)})^2 v^N-\mathscr{T}_{(d^\epsilon,P_{d^\epsilon}^*)}v^N
  )
  +\cdots$, applying the contraction property of $\mathscr{T}_{(d^\epsilon,P_{d^\epsilon}^*)}$ yields
  \begin{equation*}
  \begin{split}
  &\|(\mathscr{T}_{(d^\epsilon,P_{d^\epsilon}^*)})^\infty v^N- \mathscr{T}_{(d^\epsilon,P_{d^\epsilon}^*)}v^N\|
  \leq \lambda \|(\mathscr{T}_{(d^\epsilon,P_{d^\epsilon}^*)})^\infty v^N- v^N\|
  \\&
  \leq
  \lambda\sum_{\imath=0}^{\infty} \lambda^{\imath}\|\mathscr{T}_{(d^\epsilon,P_{d^\epsilon}^*)}v^N-v^N \|<\epsilon/2.
  \end{split}
  \end{equation*}
  Note from~(\ref{eqn23}) and the definition of $d^\epsilon$ and $P_{d^\epsilon}^*$ that the value function $v^\epsilon:=v_{(d^\epsilon,P_{d^\epsilon}^*)}$ corresponding to $d^\epsilon$ under the worst-case transition probability matrix $P_{d^\epsilon}^*$ can be given $v^\epsilon=(\mathscr{T}_{(d^\epsilon,P_{d^\epsilon}^*)})^\infty v^N$.
%
%
Therefore, $\|v^\epsilon -v^*\|\leq \|(\mathscr{T}_{(d^\epsilon,P_{d^\epsilon}^*)})^\infty v^N - \mathscr{T}_{(d^\epsilon,P_{d^\epsilon}^*)}v^N \|+\|\mathscr{T}_{(d^\epsilon,P_{d^\epsilon}^*)}v^N-v^*\|<\epsilon$.
\par
Finally, we derive the result shown in part~(c). From Proposition~\ref{prop1}, one can first get that there exists a $P'_{d_{v^*}}\in \mathcal{P}_{d_{v^*}}$ for the $v^*$-improving decision rule $d_{v^*}\in \mathcal{D}_{v^*}$ such that
  $
  \mathscr{B}v^t\geq \mathscr{B}v^* +({Q'_{d_{v^*}}}^{-1}R'_{d_{v^*}}-I)(v^t-v^*)
  $,
   where $(Q'_{d_{v^*}},R'_{d_{v^*}})$ is the GS regular splitting of $I-\lambda P'_{d_{v^*}}$.
   By applying this inequality and noticing $\mathscr{B}v^*=0$, it follows from (\ref{eqn22}) and the definition of $\mathscr{W}^{M_t}$ in~(\ref{eq22n}) that
  \begin{equation*}
  \begin{split}
  &0\leq v^* - v^{t+1}
  = v^* - \mathscr{W}^{M_t}v^{t}
  \\& =
  v^* - v^t - \sum_{\varsigma=0}^{M_t}\left({Q_{d_{v^t}}^*}^{-1} R_{d_{v^t}}^*\right)^{\varsigma}\left( \mathscr{B}v^t  \right)
  \leq
  v^* - v^t
  \\& +\sum_{\varsigma=0}^{M_t}\left({Q_{d_{v^t}}^*}^{-1} R_{d_{v^t}}^*\right)^{\varsigma}
    \left({Q'_{d_{v^*}}}^{-1} R'_{d_{v^*}}-I\right)
    (v^*-v^t)
    \\& =\left( {Q_{d_{v^t}}^*}^{-1} R_{d_{v^t}}^* -
     {Q'_{d_{v^*}}}^{-1} R'_{d_{v^*}} \right)
     \sum_{\varsigma=0}^{M_t-1} \left({Q_{d_{v^t}}^*}^{-1} R_{d_{v^t}}^*\right)^{\varsigma}
     \\&
     \times(v^t-v^*)
     -\left({Q'_{d_{v^*}}}^{-1}
     R'_{d_{v^*}}\right)
     \left({Q_{d_{v^t}}^*}^{-1} R_{d_{v^t}}^*\right)^{M_t}
     (v^t-v^*).
     \end{split}
     \end{equation*}
  Taking norms on both sides leads to
\begin{equation*}
\begin{split}
     &\|v^* - v^{t+1}\|\leq
     \| {Q_{d_{v^t}}^*}^{-1} R_{d_{v^t}}^* -
     {Q'_{d_{v^*}}}^{-1} R'_{d_{v^*}} \|\frac{1-\lambda^{M_t}}{1-\lambda}
     \\&
     \times
     \| v^t - v^* \|
     + \lambda^{M_t+1} \| v^t - v^*  \|
     =\| v^t - v^*  \|
     \\&
     \times
     \left(
     \|{Q_{d_{v^t}}^*}^{-1} R_{d_{v^t}}^* -
     {Q'_{d_{v^*}}}^{-1} R'_{d_{v^*}} \|\frac{1-\lambda^{M_t}}{1-\lambda} + \lambda^{M_t+1}
     \right),
\end{split}
\end{equation*}
{\color{red}
where $\|{Q_{d_{v^t}}^*}^{-1} R_{d_{v^t}}^* \|\leq \lambda$ and $\|{Q'_{d_{v^*}}}^{-1} R'_{d_{v^*}} \|\leq \lambda$ follow from Lemma~\ref{lm2}.
}
  In particular, if
  $\lim_{t\rightarrow \infty}
        \| P_{d_{v^t}}^{*} - P'_{d_{v^*}} \|=0$, then $\lim_{t\rightarrow \infty}     \| {Q_{d_{v^t}}^*}^{-1} R_{d_{v^t}}^* -
     {Q'_{d_{v^*}}}^{-1} R'_{d_{v^*}}  \|=0$ in view of the definition of the GS regular splitting.
     As such, there exists a $K\in \mathbb{N}$ for any $\epsilon >0$ such that
     $
     \|v^{t+1} - v^* \| \leq
     \left(\epsilon + \lambda^{M_t+1}\right)\| v^t - v^*  \|
     $
     for any $t\geq K$.
\end{IEEEproof}
{\color{red}
From this theorem, one can see that the algorithm of raTPI for $\delta=0$ can be guaranteed to
find a robust team-optimal policy at a convergence rate of near $\lambda^{M_t+1}$. Compared with the convergence rate $\lambda$ of raTVI in the degenerated case (see Theorem~\ref{thm1}), the improvement is in exponential order. Using a similar argument to that in the degenerated case, we next extend the results for $\delta=0$ to
the inexact case where $ 0\leq \delta<(1-\lambda)^2\epsilon/2\lambda(1+\lambda)$.
}
To distinguish from the notations used in the exact case $\delta=0$, we denote the sequence generated by Algorithm~1 in the inexact case by $\{\tilde{v}^t\}$.  Moreover, for those intermediate variables given in~(\ref{eqA2}) and (\ref{eqA4}), we use $\tilde{u}_0^t(s)$ and $\tilde{u}^t_{\varsigma+1}(s)$ to represent $\tilde{u}_0^t(s)=\max_{a\in\mathcal{A}}\tilde{\rho}_{(s,a)}
(\tilde{v}^t)$
and $\tilde{u}^t_{\varsigma+1}(s)=\tilde{\rho}_{(s,\tilde{d}_{t+1}(s))}
(\tilde{u}^t_\varsigma)$ in the inexact case, respectively, where $\tilde{d}_{t+1}(s)\in\arg\max_{a\in \mathcal{A}} \tilde{\rho}_{(s,a)}(\tilde{v}^t)$ for $\forall s\in \mathcal{S}$. In contrast,
we still use those notations $\{v^t\}$, $u_0^t(s)$,
$u^t_{\varsigma+1}(s)$, and $d_{t+1}(s)$ for $\forall s\in \mathcal{S}$ given in Algorithm~1 to represent the corresponding notions in the exact case. As such, we present the convergence result of raTPI in its general form by the following theorem.
\begin{theorem}\label{thm4}
For any non-negative integer sequence $\{M_t\}_{t\in \mathbb{N}}$, if the tolerance of approximation calculations is within $\lambda\delta$ for $\delta\in [0,(1-\lambda)^2\epsilon/2\lambda(1+\lambda))$ in Algorithm~1, then for any $v^0\in \mathcal{V}_{\mathscr{B}}$,
(a)~if $\tilde{v}^0=v^0$,
    then $v^t-\theta\mathbf{1}\leq \tilde{v}^t \leq v^t+\theta\mathbf{1}$ and
    $u^t_{\varsigma}-\theta\mathbf{1}\leq \tilde{u}^t_{\varsigma} \leq u^t_{\varsigma}+\theta\mathbf{1}, \quad \varsigma=0,1,\ldots,M_t$
    hold for $\forall t\in \mathbb{N}$, where $\theta:={\lambda \delta}/(1-\lambda)\in [0,\epsilon/2)$; and (b)~those results shown in the part~(b) of Theorem~\ref{thm1} still hold.
\end{theorem}

{\color{red}
The proof of this theorem is similar to those of Lemma~\ref{lm3} and Theorem~\ref{thm2}. See Appendix~\ref{app_G} for more details. From this theorem, we have demonstrated our final result that the algorithm of raTPI in its general form can be guaranteed to converge to the robust team-optimal value function, and accordingly can terminate within a finite number of iterations with an $\epsilon$-robust team-optimal policy. Finally, we give a remark to elaborate on the choice of initial iterative value functions and approximation methods.
}

\begin{remark}
The initialization condition $v^0\in \mathcal{V}_{\mathscr{B}}$ can be satisfied easily in practice. For example, selecting $v^0_s\leq (1-\lambda)^{-1}\min_{s^k,s^l\in \mathcal{S},a\in \mathcal{A}}r(s^k,a,s^l)$ for $\forall s\in \mathcal{S}$ will ensure that $v^0\in \mathcal{V}_{\mathscr{B}}$. Moreover, to implement the approximation computations in Algorithm~1, many state-of-the-art techniques can be used, such as deep neural networks or empirical samples~\cite{mannor2019data,scherrer2015approximate}. In particular, when the number of samples is large enough, it has been proven that a good enough approximate solution can be obtained by sampling-based methods~\cite{munos2008finite}.
\end{remark}
\par
To validate these theoretical results, we present some numerical simulations in the following section.
\section{Simulations}
To numerically demonstrate the effectiveness of Algorithm~1, we here generalize the game model of sequential social dilemmas in~\cite{huang2020learning} to the scenario of incomplete information, and we refer to it as \emph{Robust Sequential Social Dilemmas} (RSSDs). Formally, we consider an $n$-player Markov game where every player can only select one of the two actions, cooperation (C) and defection (D), from the action set $\mathcal{A}=\{\text{C},\text{D}\}$.
{\color{red}In response to players' actions,}
the game system will transit from the current state to a new one at the next time. More specifically, if there are $\hbar \in \{0,1,2,\ldots,n\}$ players choosing action C in the state $s^k\in \mathcal{S}$, the state at the next time will change to $s^l\in \mathcal{S}$ with probability $p(s^l\mid s^k,\hbar)$, where $p(s^l\mid s^k,\hbar)$
{\color{red}
is not predetermined but rather lies in a discrete and finite uncertain set} $\mathcal{P}_{\hbar}(l\mid k)$. {\color{red}As the consequence of players' actions and state transitions,} those players who choose action C (resp. D) will get a bounded payoff $\mathfrak{a}_\hbar(s^l\mid s^k)\in \mathbb{R}$  (resp. $\mathfrak{b}_\hbar(s^l\mid s^k)\in \mathbb{R}$) as a function of $\hbar$, $s^k$, and $s^l$. To adhere to the existence of dilemmas, we assume as in the canonical multi-player social dilemmas~\cite{Govaert} that (i) $\mathfrak{a}_{\hbar+1}(s^l\mid s^k)\geq \mathfrak{a}_{\hbar}(s^l\mid s^k)$ and $\mathfrak{b}_{\hbar+1}(s^l\mid s^k)\geq \mathfrak{b}_{\hbar}(s^l\mid s^k)$;
(ii) $\mathfrak{b}_{\hbar}(s^l\mid s^k)>\mathfrak{a}_{\hbar}(s^l\mid s^k)$; and
(iii) $\mathfrak{a}_{n}(s^l\mid s^k)>\mathfrak{b}_{0}(s^l\mid s^k)$ for all $\hbar$ and $\forall s^k,s^l\in \mathcal{S}$. Condition (i) states that players' payoffs increase with the number of C players in the group, whereas condition (ii) implies that within any mixed group, those C players always have a strictly lower payoff than that of D players. These two conditions indicate that taking action C is an altruistic behavior because C players entail a potential cost to improve other agents' benefits. In contrast, condition (iii) shows that mutual cooperation is more beneficial than mutual defection. Hence, for maximizing the gains of the whole group, all players should uniformly choose C. However, under the hypothesis of Homo economicus, each rational player will be tempted by myopic interests to take action D, thereby leading to the existence of social dilemmas.
By extending social dilemmas to the scenario of incomplete information, the current model incorporates the prototypical multi-player social dilemmas~\cite{Govaert} as an extreme case where the set of states is a singleton. Also, since the uncertainty of transition probabilities implies that there are some parameters of games unknown to players, it generalizes social dilemmas to the scenario of robust games~\cite{aghassi2006robust,kardecs2011discounted}
\par
\begin{figure}
  \centering
  \includegraphics[width=\hsize]{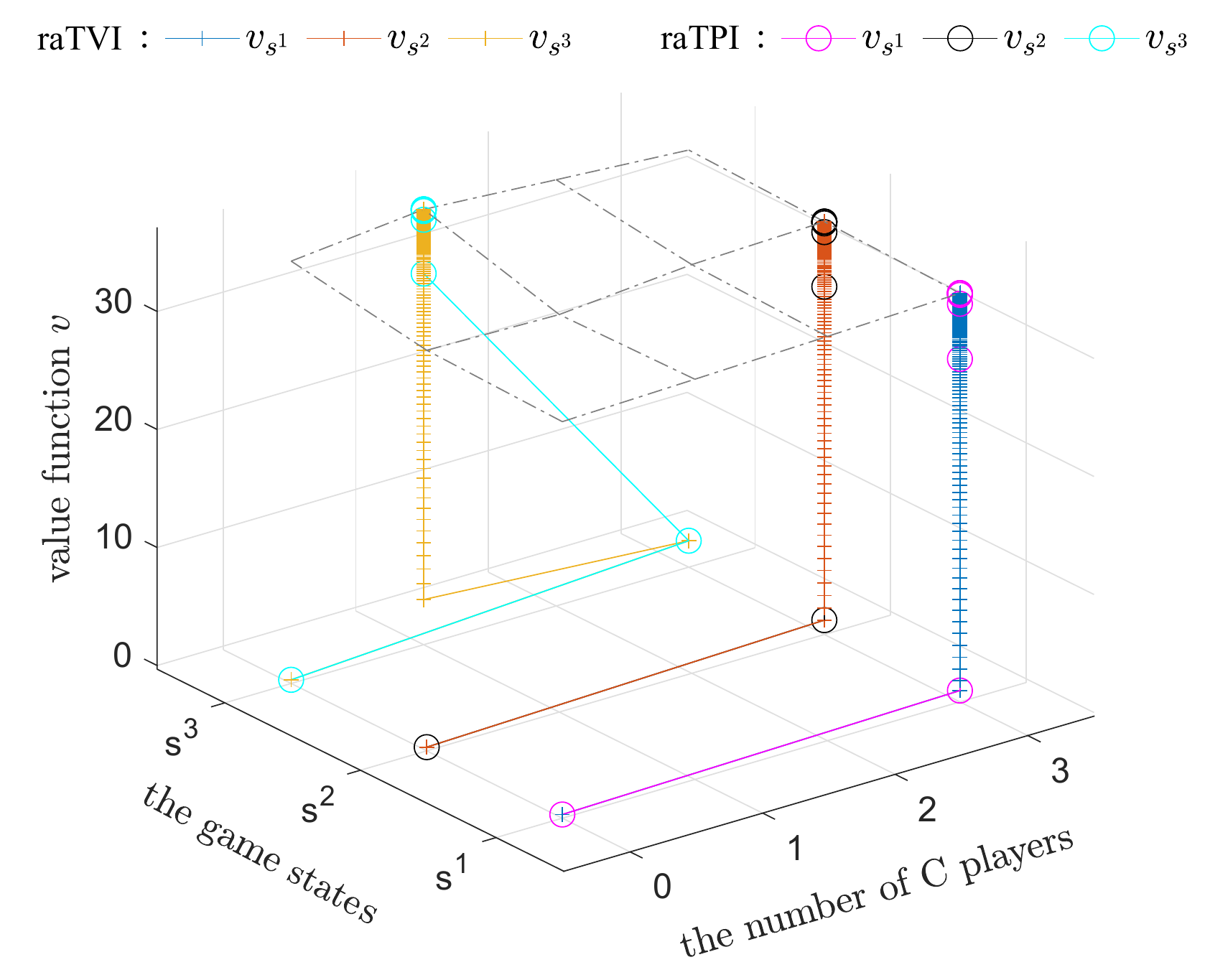}\\
  \caption{The iterative process of the value functions generated by raTVI and raTPI. {\color{red}Each point on lines is drawn by the iterative data.} The dashed lattice illustrates the values of $\rho_{(s,a)}(v^t)$ for $\forall (s,a)\in \mathcal{S}\times\mathcal{A}$ at the termination time. } \label{fig1}
\end{figure}

{\color{red}
In contrast to the underlying non-cooperative setting in conventional social dilemmas, we here assume that although every player has its individual payoff, all players aim to maximize the long-term benefit of the whole group. Such an assumption is to some extent reasonable in some social dilemmas encountered by humans, such as the coalition of nations when facing global warming. Under this assumption, the game model of RSSDs then becomes a specific instance of robust team Markov games. Therefore, we can apply Algorithm~1 to seek the robust team-optimal policy. Specifically, we consider a game of RSSDs with the state set of three elements, $\mathcal{S}=\{s^1, s^2, s^3\}$.
}
In state $s^1$, players play a public goods game and players' payoffs are calculated by $\mathfrak{a}_{\hbar}(s^l\mid s^1)=\hbar r_{s^l}c/n-c$ and $\mathfrak{b}_{\hbar}(s^l\mid s^1)=\hbar r_{s^l}c/n$ for $\forall s^l\in \mathcal{S}$, where $c$ is the cost of cooperation and $r_{s^l}\in(c,n)$ is the synergy factor dependent on the state $s^l$ at the next time. In state $s^2$, an $n$-player stag-hunt game is played, where players' payoffs are calculated by $\mathfrak{a}_{\hbar}(s^l\mid s^2)=\hbar r_{s^l}c/n-c$ and $\mathfrak{b}_{\hbar}(s^l\mid s^2)=\hbar r_{s^l}c/n$ if $\hbar$ is no less than the threshold $Z$; and otherwise $\mathfrak{a}_{\hbar}(s^l\mid s^2)=-c$ and $\mathfrak{b}_{\hbar}(s^l\mid s^2)=0$. In state $s^3$, an $n$-player snow-drift game is played, and players' payoffs are calculated by $\mathfrak{a}_{\hbar}(s^l\mid s^3)=\vartheta_{s^l}-c/\hbar$ and $\mathfrak{b}_{\hbar}(s^l\mid s^3)=\vartheta_{s^l}$ if $\hbar>0$, and otherwise $\mathfrak{a}_{\hbar}(s^l\mid s^3)=\mathfrak{b}_{\hbar}(s^l\mid s^3)=0$, where $\vartheta_{s^l}$, as a function of the state $s^l$ at the next time, is the benefit of players when there exists at least one player choosing action C in the group. (see~\cite{huang2020learning} and references therein for more details on these three games.) The simulation results are shown in~Fig.~\ref{fig1} and Table~\ref{tab1}, where the model parameters of games are given by $r_{s^1}=1.5$, $r_{s^2}=1.8$, $r_{s^3}=2.2$, $\vartheta_{s^l}=r_{s^l}$ for $\forall s^l\in \mathcal{S}$, $n=3$, $c=1.0$, and $\mathcal{P}_{\hbar}(l\mid k)=\{1-\mu \hbar\}_{\mu\in\{0.1,0.2,0.3\}}$ if $k=l$ and otherwise $\mathcal{P}_{\hbar}(l\mid k)=\{\mu \hbar/2\}_{\mu\in\{0.1,0.2,0.3\}}$ for $\forall \hbar\in\{0,\ldots,n\}$. Moreover, the algorithm parameters are $\epsilon=10^{-5}$ and $\delta\approx(1-\lambda)^2\epsilon/2\lambda(1+\lambda)$, where a fixed $\lambda=0.97$ is adopted in Fig.~\ref{fig1} while different $\lambda$ values are used in Table~\ref{tab1}. From Fig.~\ref{fig1} and Table~\ref{tab1}, one can see that both raTVI and raTPI are able to effectively find the robust team-optimal policy, and their convergence rates are faster than rVI and rMPI, respectively.
{\color{red}
Although the computational model adopted here is minimal, the improvement is noticeable. In problems with large action and state sets, moreover, the performance of Algorithm~1 will be in general more competitive.
}

\begin{table}[htbp]
  \centering
  \def~{\hphantom{0}}
  \caption{The number of iterations for different algorithms to find the $\epsilon$-robust team-optimal policy.}
  \label{tab1}
  \begin{tabular*}{\hsize}{@{\extracolsep{\fill}}cccccc}
    \toprule
    \multirow{2}{*}{Algorithms}&
    \multicolumn{5}{c}{the magnitude of $\lambda$}\cr
     \cmidrule{2-6}
    &0.95&0.96&0.97&0.98&0.99\cr
    \midrule
    rVI&298&380&519&802&1679\cr
    raTVI&\textbf{258}&\textbf{328}&\textbf{446}&\textbf{690}
    &\textbf{1442}\cr
    rMPI&7&9&12&17&34\cr
    raTPI&\textbf{7}&\textbf{8}&\textbf{10}&\textbf{15}
    &\textbf{30}\cr
    \bottomrule
  \end{tabular*}
\end{table}

%

\section{Conclusions and future work}
{\color{red}
In this work, we have addressed a sequential decision-making problem of a cooperative team in stochastic uncertain environments,} and accordingly proposed a robust version of team Markov games by relaxing the complete information assumption, in which players do not know the accurate transition probabilities of states, but rather are commonly aware of an uncertainty set. Without assuming a prior probability distribution over the uncertainty set, we consider that players adopt robust optimization methods to update their strategies. To characterize the optimality of team decisions, we have also proposed a solution concept named robust team-optimal
policy, and meanwhile developed a robust iterative learning algorithm to seek it. Under mild conditions, we have presented the convergence analysis of the algorithm and calculated its convergence rates. {\color{red}By numerical simulations, moreover, we have demonstrated the effectiveness of these theoretical results in a game of RSSDs.
However, the paper has left some open questions, which point to future directions: how to achieve decentralized or individual robust learning with theoretical guarantees, especially in network/graphic games; and how to achieve sample efficiency by data-driven robust learning or design robust policy search algorithms by function approximation. Also, we will be interested in applying RSSDs to study real-world issues, such as the evolution of social power~\cite{tian2021social} and human decisions in social diffusion~\cite{ye2021collective}, in the future.}

\appendix
{\color{red}
In this section, we provide some proof details of the results shown in the main text.
}
\subsection{Proof of Lemma~1}
\label{app_A}
{\color{red}
We establish part~(a) by first proving that $\mathscr{Y}$ is a contraction mapping on $\mathcal{V}$.} For any given $u,v\in \mathcal{V}$, we start with considering those states $s \in \mathcal{S}$ satisfying $(\mathscr{Y} v-\mathscr{Y} u)(s)\geq 0$, where $(\mathscr{Y} v-\mathscr{Y} u)(s)$ is the component of $\mathscr{Y} v-\mathscr{Y} u$ corresponding to $s\in \mathcal{S}$.
Select $
d_v\in \arg\max_{d\in \mathcal{D}}\left\{
\min_{P_d\in \mathcal{P}_d}\left[Q_d^{-1}r_{(d,P_d)}+
Q_d^{-1}R_d v
\right]
\right\}
$.
{\color{red}
It then follows from the definition of $\mathscr{Y}$ in~(\ref{eq14}) that
}
\begin{equation}\label{eqA0}
\begin{split}
& (\mathscr{Y} v- \mathscr{Y} u)(s)
 \leq
 \left\{\min_{P_{d_v}\in \mathcal{P}_{d_v}}\left\{
Q_{d_v}^{-1}r_{({d_v},P_{d_v})}
+ Q_{d_v}^{-1}R_{d_v} v
\right\}\right. \\
&
\left.-\min_{P_{d_v}\in \mathcal{P}_{d_v}}\left\{
Q_{d_v}^{-1}r_{({d_v},P_{d_v})} + Q_{d_v}^{-1}R_{d_v} u
\right\}
\right\}(s).
\end{split}
\end{equation}
Let $
\hat{P}_{d_v}\in \arg\min_{P_{d_v}\in \mathcal{P}_{d_v}}
\left\{
Q_{d_v}^{-1}r_{({d_v},P_{d_v})} + Q_{d_v}^{-1}R_{d_v} u
\right\}
$.
{\color{red}
Then, the right-hand side of~(\ref{eqA0}) will be no larger than
}
\begin{equation*}
\begin{split}
&\left\{\left(
\hat{Q}_{d_v}^{-1}r_{({d_v},\hat{P}_{d_v})} + \hat{Q}_{d_v}^{-1}\hat{R}_{d_v} v
\right)
-
\left(
\hat{Q}_{d_v}^{-1}r_{({d_v},\hat{P}_{d_v})} \right.\right.
\\
&
\left.\left.
+ \hat{Q}_{d_v}^{-1}\hat{R}_{d_v} u
\right)\right\}(s)
=
\left\{\hat{Q}_{d_v}^{-1}\hat{R}_{d_v} (v-u)\right\}(s)
\\
&
\leq \|\hat{Q}_{d_v}^{-1}\hat{R}_{d_v} (v-u)\|
\leq \|\hat{Q}_{d_v}^{-1}\hat{R}_{d_v}\| \|v-u \|
  \leq \alpha \|v-u \|,
\end{split}
\end{equation*}
where
$(\hat{Q}_{d_v},\hat{R}_{d_v})$ is the corresponding regular splitting of $I-\lambda \hat{P}_{d_v}$. Consequently, $0 \leq (\mathscr{Y} v- \mathscr{Y} u)(s)\leq \alpha\|v-u\|$.
Similarly, for those states $s \in \mathcal{S}$ satisfying $ (\mathscr{Y} v-\mathscr{Y} u)(s)\leq 0$, one can obtain $0 \leq (\mathscr{Y} u- \mathscr{Y} v)(s)\leq \alpha\|u-v\|$ by the same argument. Therefore, $\|\mathscr{Y} v- \mathscr{Y} u\|\leq \alpha\|v-u\|$. Since $\alpha<1$, $\mathscr{Y}$ is a contraction mapping on $\mathcal{V}$. From the Banach fixed-point theorem (see Theorem~$6.2.3$~in~\cite{puterman2014markov}), it follows that $\{v^t\}$ converges in norm to the unique fixed point $\tilde{v}^{*}$ of $\mathscr{Y}$.
\par
We next show that $\tilde{v}^{*}=v^{*}$. Note that $\tilde{v}^{*}$ is the unique fixed point of $\mathscr{Y}$. Then, arbitrarily given a $d\in \mathcal{D}$,  we have
$
\tilde{v}^*=\mathscr{Y} \tilde{v}^* \geq \min_{P_d\in \mathcal{P}_d}\left\{Q_d^{-1}r_{(d,P_d)} + Q_d^{-1}R_d \tilde{v}^* \right\}
$ from the definition of $\mathscr{Y}$ in~(\ref{eq14}).
Also, for any $\epsilon>0$, there exists a $P_d^{\epsilon}\in \mathcal{P}_d$ such that
$
\min_{P_d\in \mathcal{P}_d}\left\{Q_d^{-1}r_{(d,P_d)} + Q_d^{-1}R_d \tilde{v}^*\right\} \geq
{Q_d^{\epsilon}}^{-1}r_{(d,P_d^{\epsilon})}+
{Q_d^{\epsilon}}^{-1}R_d^{\epsilon} \tilde{v}^* - \epsilon \mathbf{1}
$,
where $(Q_d^{\epsilon},R_d^{\epsilon})$ is the corresponding regular splitting of $I-\lambda P_d^{\epsilon}$. As a result,
$
\tilde{v}^*\geq {Q_d^{\epsilon}}^{-1}r_{(d,P_d^{\epsilon})}+
{Q_d^{\epsilon}}^{-1}R_d^{\epsilon} \tilde{v}^* - \epsilon \mathbf{1}
$.
By re-arrangement, it leads to
\begin{equation}\label{Aa1}
\begin{split}
\tilde{v}^* &\geq (I-{Q_d^{\epsilon}}^{-1}R_d^{\epsilon})^{-1}
{Q_d^{\epsilon}}^{-1} r_{(d,P_d^{\epsilon})}-\epsilon (I-{Q_d^{\epsilon}}^{-1}R_d^{\epsilon})^{-1} \mathbf{1} \\
&=(Q_d^{\epsilon}-R_d^{\epsilon})^{-1} r_{(d,P_d^{\epsilon})}
-\epsilon (I-{Q_d^{\epsilon}}^{-1}R_d^{\epsilon})^{-1} \mathbf{1} \\
&= (I-\lambda P_d^{\epsilon})^{-1}r_{(d,P_d^{\epsilon})} -\epsilon \sum_{\imath=0}^{\infty} ({Q_d^{\epsilon}}^{-1}R_d^{\epsilon})^\imath \mathbf{1} \\
&= \sum_{\imath=0}^{\infty} (\lambda P_d^{\epsilon})^\imath r_{(d,P_d^{\epsilon})} -\epsilon \sum_{\imath=0}^{\infty} ({Q_d^{\epsilon}}^{-1}R_d^{\epsilon})^\imath \mathbf{1} \\
&= v_{(d,P_d^{\epsilon})} -\epsilon \sum_{\imath=0}^{\infty} ({Q_d^{\epsilon}}^{-1}R_d^{\epsilon})^\imath \mathbf{1},
\end{split}
\end{equation}
where $(I-\lambda P_d^{\epsilon})^{-1}=\sum_{\imath=0}^{\infty} (\lambda P_d^{\epsilon})^\imath $ and $(I-{Q_d^{\epsilon}}^{-1}R_d^{\epsilon})^{-1}=\sum_{\imath=0}^{\infty} ({Q_d^{\epsilon}}^{-1}R_d^{\epsilon})^\imath$ hold because $\sigma(\lambda P_d^{\epsilon})\leq \|\lambda P_d^{\epsilon}\|=\lambda<1$ and $\sigma({Q_d^{\epsilon}}^{-1}R_d^{\epsilon})\leq \|{Q_d^{\epsilon}}^{-1}R_d^{\epsilon}\|\leq \alpha <1$;
{\color{red}and the last equality follows from the first equality of~(\ref{eq5}).
}
Since for any non-negative integer $\imath$,
\begin{equation}\label{eqAa1}
\begin{split}
&({Q_d^{\epsilon}}^{-1}R_d^{\epsilon})^{\imath}\mathbf{1}\leq \mathbf{1}\sup_{k\in \{1,2,\ldots,m\}}\sum_{l=1}^{m}
|({Q_d^{\epsilon}}^{-1}R_d^{\epsilon})^{\imath}(l|k)|
\\
&
=
\mathbf{1}\|({Q_d^{\epsilon}}^{-1}R_d^{\epsilon})^{\imath}\|\leq
\mathbf{1}\|{Q_d^{\epsilon}}^{-1}R_d^{\epsilon}\|^{\imath}\leq
\alpha^{\imath}\mathbf{1},
\end{split}
\end{equation}
the inequality (\ref{Aa1}) can be further given by
$
\tilde{v}^* \geq v_{(d,P_d^{\epsilon})} -\epsilon \sum_{\imath=0}^{\infty} ({Q_d^{\epsilon}}^{-1}R_d^{\epsilon})^\imath \mathbf{1} \geq
v_{(d,P_d^{\epsilon})} -\epsilon \sum_{\imath=0}^{\infty}\alpha^{\imath}\mathbf{1} =
v_{(d,P_d^{\epsilon})} -{\epsilon}\mathbf{1}/(1-\alpha).
$
As $\epsilon>0$ is arbitrary, $\tilde{v}^*\geq v_{(d,P_d^{\epsilon})} \geq \min_{P_d\in\mathcal{P}_d} v_{(d,P_d)}$. In addition, note that $d\in \mathcal{D}$ is arbitrary. Thus, $\tilde{v}^*\geq v^*=\max_{d\in \mathcal{D}} \min_{P_d\in \mathcal{P}_d} v_{(d,P_d)}$.
\par
On the other hand, since from the definition of $\mathscr{Y}$, $\tilde{v}^*=\mathscr{Y} \tilde{v}^*=\max_{d\in \mathcal{D}}\min_{P_d\in\mathcal{P}_d} \left\{Q_d^{-1}r_{(d,P_d)}+Q_d^{-1}R_d \tilde{v}^*\right\}$, there exists a $\tilde{d}\in \mathcal{D}$ for any $\epsilon>0$ such that $\tilde{v}^*\leq \min_{P_{\tilde{d}}\in\mathcal{P}_{\tilde{d}}} \left\{ Q_{\tilde{d}}^{-1}r_{(\tilde{d},P_{\tilde{d}})}+Q_{\tilde{d}}^{-1}R_{\tilde{d}} \tilde{v}^*\right\} +\epsilon \mathbf{1}$. Thus, for any $\tilde{P}_{\tilde{d}} \in \mathcal{P}_{\tilde{d}}$,
\begin{equation}\label{Aa2}
\tilde{v}^*\leq {\tilde{Q}}_{\tilde{d}}^{-1}r_{(\tilde{d},{\tilde{P}}_{\tilde{d}})}+{\tilde{Q}}_{\tilde{d}}^{-1}
{\tilde{R}}_{\tilde{d}} \tilde{v}^* +\epsilon \mathbf{1},
\end{equation}
where $(\tilde{Q}_{\tilde{d}},\tilde{R}_{\tilde{d}})$ is the corresponding regular splitting of $I-\lambda \tilde{P}_{\tilde{d}}$. Using a similar argument to~(\ref{Aa1}) and~(\ref{eqAa1}), one can then get $\tilde{v}^*\leq v_{(\tilde{d},\tilde{P}_{\tilde{d}})} +{\epsilon}\mathbf{1}/(1-\alpha)$ from (\ref{Aa2}). Since both $\epsilon>0$ and $\tilde{P}_{\tilde{d}} \in \mathcal{P}_{\tilde{d}}$ are arbitrary for the given $\tilde{d}\in \mathcal{D}$, $\tilde{v}^*\leq \min_{P_{\tilde{d}}\in \mathcal{P}_{\tilde{d}}}v_{(\tilde{d},P_{\tilde{d}})}$. As such, $\tilde{v}^*\leq \max_{d\in \mathcal{D}} \min_{P_d\in \mathcal{P}_d} v_{(d,P_d)} =v^*$. Therefore, we have $\tilde{v}^{*}=v^{*}$.
\par
Finally, we establish part (b). Note first from part~(a) that $\mathscr{Y}$ is a contraction mapping on $\mathcal{V}$ with constant $\alpha$ and its unique fixed point is $v^*$. Thus,
$
   \|v^{t+1}-v^*\|= \|\mathscr{Y} v^t - \mathscr{Y} v^*\|\leq \alpha\|v^t - v^*\|$ for $\forall t\in \mathbb{N}
$.
Then, using this inequality recursively yields $\|v^t-v^*\|\leq \alpha^t \|v^0 - v^*\|$. Dividing both sides by $\|v^0 - v^*\|$ and taking the $t$-th root show that
$
\beta:=\limsup_{t\rightarrow\infty}\left[
      {\|v^t-v^*\|}/{\|v^0-v^*\|}
      \right]^{1/t}\leq \alpha.
$
From this inequality,
$
\limsup_{t\rightarrow\infty}{\|v^t-v^*\|}/{\beta^t}\leq
\|v^0-v^*\|<\infty
$
can be immediately obtained.
{\color{red}
From Definition~\ref{def3}, the proof of
part (b) is completed.
}

\subsection{Proof of Lemma~3}
\label{app_B}
{\color{red}
We establish the proof by induction.} For $t=0$, (\ref{eqn14}) follows from the assumption $v^0=\tilde{v}^0$. Suppose now that (\ref{eqn14}) is satisfied for $t=\kappa\geq 0$, and we derive the result for $t=\kappa+1$.
{\color{red}
Note first that
for any given $a\in \mathcal{A}$ in state $s^k\in\mathcal{S}$, the term inside the curly brace in (\ref{eqn9}) is monotonically increasing with respect to $v^t$ for any given $p(\cdot|s^k,a)\in \mathcal{P}_a(\cdot|k)$. Moreover, since the minimization operation in (\ref{eqn9}) over the probability distribution will not change the monotonicity, $\rho_{(s^k,a)}(v^t)$ is monotonically increasing with respect to $v^t$.
}
For any given $a\in \mathcal{A}$ in state $s^1$, substituting each term of (\ref{eqn14}) into $\rho_{(s^1,a)}(\cdot)$ then yields
\begin{equation*}
\begin{split}
\rho_{(s^1,a)}(v^\kappa-\theta\mathbf{1}) &= \rho_{(s^1,a)}(v^\kappa)-\lambda \theta
\leq \rho_{(s^1,a)}(\tilde{v}^\kappa)
\\
&
\leq \rho_{(s^1,a)}(v^\kappa)+\lambda \theta
=\rho_{(s^1,a)}(v^\kappa+\theta\mathbf{1}).
\end{split}
\end{equation*}
Taking the maximum of each term of this inequality over $a\in \mathcal{A}$, one can get
\begin{equation}\label{eqn16}
    v_{s^1}^{\kappa+1}-\lambda \theta\leq \max_{a\in \mathcal{A}}
    \rho_{(s^1,a)}(\tilde{v}^\kappa) \leq
    v_{s^1}^{\kappa+1}+\lambda \theta,
\end{equation}
{\color{red}where $v_{s^1}^{\kappa+1}=\max_{a\in \mathcal{A}}\rho_{(s^1,a)}(v^\kappa)$ follows from~(\ref{eqn11}).
}
Since from~(\ref{eqA1}), $\tilde{\rho}_{(s^1,a)}(\tilde{v}^\kappa)-\lambda \delta\leq \rho_{(s^1,a)}(\tilde{v}^\kappa)\leq \tilde{\rho}_{(s^1,a)}(\tilde{v}^\kappa)+\lambda \delta$ for $\tilde{v}^\kappa$, taking the maximum of each term subject to $a\in \mathcal{A}$ yields
\begin{equation}\label{eqn17}
  \tilde{v}^{\kappa+1}_{s^1}-\lambda \delta \leq \max_{a\in \mathcal{A}}\rho_{(s^1,a)}(\tilde{v}^\kappa) \leq
  \tilde{v}^{\kappa+1}_{s^1}+\lambda \delta,
\end{equation}
{\color{red}
where
$\tilde{v}_{s^1}^{\kappa+1}=\max_{a\in \mathcal{A}}\tilde{\rho}_{(s^1,a)}(\tilde{v}^\kappa)$ follows from~(\ref{eqA2}) and the step~3(f) of Algorithm~1 for $\delta\in [0,(1-\lambda)^2\epsilon/2\lambda(1+\lambda))$ in the case of $M_t=0$.
It then follows from~(\ref{eqn16}) and (\ref{eqn17}) that
}
\begin{equation}\label{eqn18}
v^{\kappa+1}_{s^1}-\lambda(\delta+\theta)\leq
\tilde{v}^{\kappa+1}_{s^1}\leq v^{\kappa+1}_{s^1}+\lambda(\delta+\theta).
\end{equation}
Since $\lambda(\delta+\theta)=\theta$, (\ref{eqn18}) is changed to
$
v^{\kappa+1}_{s^1}-\theta\leq
\tilde{v}^{\kappa+1}_{s^1}\leq v^{\kappa+1}_{s^1}+\theta.
$
Combing this inequality with the hypothesis $v^\kappa_{s^j}-\theta\leq \tilde{v}^\kappa_{s^j} \leq v^\kappa_{s^j}+\theta$, $j=2,3,\ldots,m$, and from~(\ref{eqn9}), one can further obtain
$
\rho_{(s^2,a)}(v^\kappa)-\lambda \theta\leq \rho_{(s^2,a)}(\tilde{v}^\kappa)\leq \rho_{(s^2,a)}(v^\kappa)+\lambda \theta
$ for $\forall a\in \mathcal{A}$ by a similar argument to the above process.
Taking the maximum of each term subject to
$a\in \mathcal{A}$ then leads to
\begin{equation}\label{eqC3}
v_{s^2}^{\kappa+1}-\lambda \theta\leq \max_{a\in \mathcal{A}} \rho_{(s^2,a)}
(\tilde{v}^\kappa) \leq v_{s^2}^{\kappa+1}+\lambda \theta.
\end{equation}
Likewise, for $s^2\in \mathcal{S}$, since $\tilde{\rho}_{(s^2,a)}(\tilde{v}^\kappa)-\lambda \delta\leq \rho_{(s^2,a)}(\tilde{v}^\kappa)\leq \tilde{\rho}_{(s^2,a)}(\tilde{v}^\kappa)+\lambda \delta$ for $\tilde{v}^\kappa$ from~(\ref{eqA1}), taking the maximum of each term over $a\in\mathcal{A}$ leads to
\begin{equation}\label{eqC4}
\tilde{v}^{\kappa+1}_{s^2}-\lambda \delta \leq \max_{a\in \mathcal{A}}\rho_{(s^2,a)}(\tilde{v}^\kappa) \leq
  \tilde{v}^{\kappa+1}_{s^2}+\lambda \delta.
\end{equation}
Using (\ref{eqC3}) and (\ref{eqC4}), one can similarly obtain
$
v^{\kappa+1}_{s^2}-\theta\leq
\tilde{v}^{\kappa+1}_{s^2}\leq v^{\kappa+1}_{s^2}+\theta
$.
Then, applying the same argument for $s^j$, $j=3,4,\ldots,m$ by recursion, one can get $v^{\kappa+1}_{s^j}-\theta \leq
{\tilde{v}}^{\kappa+1}_{s^j}\leq v^{\kappa+1}_{s^j}+\theta$. Therefore,  $v^{\kappa+1}-\theta\mathbf{1}\leq \tilde{v}^{\kappa+1} \leq v^{\kappa+1}+\theta\mathbf{1}$. The proof is completed.

\subsection{Proof of Proposition~1}
\label{app_C}
For any given $u,v\in \mathcal{V}$ and $d_v\in \mathcal{D}_v$, from the definition of $\mathscr{B}$ in~(\ref{eq18}), one can get
\begin{equation}\label{eqC1}
  \mathscr{B}u\geq \min_{P_{d_v}\in \mathcal{P}_{d_v}}\left\{
  Q_{d_v}^{-1}r_{({d_v},P_{d_v})} +
  Q_{d_v}^{-1}R_{d_v} u - u\right\}
\end{equation}
  and
\begin{equation}\label{eqC2}
  \mathscr{B}v =\min_{P_{d_v}\in \mathcal{P}_{d_v}}\left\{
  Q_{d_v}^{-1}r_{({d_v},P_{d_v})}+
  Q_{d_v}^{-1}R_{d_v} v  - v\right\}.
\end{equation}
Moreover, for any $\epsilon>0$, there exists a $P'_{d_v}\in \mathcal{P}_{d_v}$ such that the term on the right-hand side of~(\ref{eqC1}) is no less than
  $ {Q'_{d_v}}^{-1}r_{({d_v},P'_{d_v})} +
  {Q'_{d_v}}^{-1}R'_{d_v} u - u-\epsilon \mathbf{1}$, and meanwhile the term on the right-hand side of~(\ref{eqC2}) is no larger than $ {Q'_{d_v}}^{-1}r_{({d_v},P'_{d_v})} +
  {Q'_{d_v}}^{-1}R'_{d_v} v - v$,
where $(Q'_{d_v},R'_{d_v})$ is the GS regular splitting of $I-\lambda P'_{d_v}$. Then,
subtracting $\mathscr{B}v$ from $\mathscr{B}u$ yields
  $\mathscr{B}u - \mathscr{B}v \geq ({Q'_{d_v}}^{-1}R'_{d_v}-I)(u-v)-\epsilon \mathbf{1}$. Note that $\epsilon>0$ is arbitrary. Therefore, $\mathscr{B}u \geq \mathscr{B}v + ({Q'_{d_v}}^{-1}R'_{d_v}-I)(u-v)$.

\subsection{Proof of Lemma~\ref{lm4}}
\label{app_D}
Since for any given $d\in \mathcal{D}$ and $P_d\in \mathcal{P}_d$, the $(Q_d,R_d)$ given in the definition~(\ref{eq23}) of $\mathscr{U}^M$ is the GS regular splitting of $I-\lambda P_d$,  $\|{Q_d}^{-1}R_d\| \leq \| I^{-1} (\lambda P_d)\|=\lambda<1$ from Lemma~\ref{lm2}. It implies that $\sup_{d\in \mathcal{D},P_d\in \mathcal{P}_d}\|Q_d^{-1}R_d\|\leq\lambda$. Thus, we can adopt a similar argument to the proof of Lemma~1(a) to prove part~(a). For any given $u,v\in \mathcal{V}$, we begin with considering those states $s \in \mathcal{S}$ satisfying $(\mathscr{U}^M v-\mathscr{U}^M u)(s)\geq 0$. Select
  $
  d'_v\in \arg\max_{d\in \mathcal{D}}\left\{\max_{P_d\in \mathcal{P}_d} \Phi(d,P_d,v)\right\}.
  $
Then, from the definition of $\mathscr{U}^M$ in~(\ref{eq23}), we have
$$
\mathscr{U}^M v(s)=\{\max_{P_{d_v'}\in \mathcal{P}_{d_v'}}
  \Phi(d'_v,P_{d'_v},v) \}(s),
$$
and
$$
\mathscr{U}^M u(s) \geq
  \{ \max_{P_{d_v'}\in \mathcal{P}_{d_v'}}
  \Phi(d'_v,P_{d'_v},u) \}(s)
.
$$
Moreover,
since there exists a $P'_{d'_v}\in \mathcal{P}_{d'_v}$ for any $\epsilon>0$ such that
$\max_{P_{d_v'}\in \mathcal{P}_{d_v'}}\Phi(d'_v,P_{d'_v},v)\leq \Phi(d'_v,P'_{d'_v},v)+\epsilon\mathbf{1}$ and meanwhile $\max_{P_{d_v'}\in \mathcal{P}_{d_v'}}
  \Phi(d'_v,P_{d'_v},u) \geq \Phi(d'_v,P'_{d'_v},u)$, subtracting $\mathscr{U}^M u(s)$ from $\mathscr{U}^M v(s)$ leads to
  $$
  (\mathscr{U}^M v - \mathscr{U}^M u)(s) \leq \left[ \Phi(d'_v,P'_{d'_v},v)+\epsilon\mathbf{1} - \Phi(d'_v,P'_{d'_v},u) \right](s).
  $$
  Note that $\epsilon>0$ is arbitrary. Therefore, from the definition of $\Phi(\cdot)$ in~(\ref{eq23}), we further have
  \begin{equation*}
  \begin{split}
  (\mathscr{U}^M v - \mathscr{U}^M u)(s)
  &\leq \left[ \Phi(d'_v,P'_{d'_v},v) - \Phi(d'_v,P'_{d'_v},u) \right](s)
  \\
  &=
  \left[
  ({Q'_{d'_v}}^{-1} R'_{d'_v})^{M+1}
  (v-u)
  \right](s),
  \end{split}
  \end{equation*}
where $(Q'_{d'_v},R'_{d'_v})$ is the GS regular splitting of $I-\lambda P'_{d'_v}$.
Moreover, note that
\begin{equation*}
\begin{split}
&\left[({Q'_{d'_v}}^{-1} R'_{d'_v})^{M+1}(v-u)\right](s)
\leq
  \|({Q'_{d'_v}}^{-1} R'_{d'_v})^{M+1}
  (v-u) \|
\\
&
\leq
  \|{Q'_{d'_v}}^{-1} R'_{d'_v}\|^{M+1} \|v-u\|
\leq \lambda^{M+1} \|v-u\|.
\end{split}
\end{equation*}
Therefore, $0 \leq (\mathscr{U}^M v - \mathscr{U}^M u)(s) \leq \lambda^{M+1} \|v-u\|$.
\par
Applying the same argument for those states $s \in \mathcal{S}$ satisfying $ (\mathscr{U}^M v-\mathscr{U}^M u)(s)\leq 0$, one can similarly obtain $0\leq (\mathscr{U}^M u-\mathscr{U}^M v)(s) \leq \lambda^{M+1}\|u-v\|$. Consequently,
$
\|\mathscr{U}^M v - \mathscr{U}^M u\| \leq \lambda^{M+1}\|v-u\|
$.
Since $\lambda^{M+1}<1$, $\mathscr{U}^M$ is a contraction mapping on $\mathcal{V}$. It then follows from the Banach fixed-point theorem that for any $\omega^0 \in \mathcal{V}$, the sequence $\{\omega^t\}$ generated by $\mathscr{U}^M$ will converge in norm to the unique fixed point $\omega^*$ of $\mathscr{U}^M$.
\par
We next show $\omega^*=v^*$. Note first that $v^*$ is the unique fixed point of $\mathscr{Y}$ from Lemma~\ref{lm1}. Then, based on the definition of $\mathscr{Y}$ in~(\ref{eq14}), we have $v^*=\mathscr{Y} v^*=\max_{d\in \mathcal{D}}\min_{P_d\in\mathcal{P}_d}
\left\{Q_d^{-1}r_{(d,P_d)}+Q_d^{-1}R_d v^*\right\}$. Let $d_{v^*}\in \mathcal{D}_{v^*}$ be a $v^*$-improving decision rule and $P_{d_{v^*}}^*\in \mathcal{P}_{d_{v^*}}^*$ be a $d_{v^*}$-decreasing transition probability matrix. As such, $v^*=
\mathscr{T}_{(d_{v^*},P_{d_{v^*}}^*)}v^*$. {\color{red} Leveraging (\ref{eqn22}) and the definition of $\mathscr{W}^M$ in (\ref{eq22n}), one can then get
}
\begin{equation*}
\begin{split}
v^* &= (\mathscr{T}_{(d_{v^*},P_{d_{v^*}}^*)})^{M+1}
v^*
=
\mathscr{W}^M v^*
\\
&
= \sum_{\varsigma=0}^{M}({Q_{d_{v^*}}^*}^{-1} R_{d_{v^*}}^*)^{\varsigma}{Q_{d_{v^*}}^*}^{-1}
  r_{(d_{v^*},P_{d_{v^*}}^*)}
\\
&
  +({Q_{d_{v^*}}^*}^{-1} R_{d_{v^*}}^*)^{M+1}v^*
  =\Phi(d_{v^*},P_{d_{v^*}}^*,v^*),
\end{split}
\end{equation*}
where $(Q_{d_{v^*}}^*,R_{d_{v^*}}^* )$ is the GS regular splitting of $I-\lambda P_{d_{v^*}}^*$.
From the definition of $\mathscr{U}^M$ in (\ref{eq23}), it then leads to $v^*\leq \mathscr{U}^M v^*$.
Applying this inequality recursively for $\imath \in \mathbb{N}$ times yields $v^*\leq (\mathscr{U}^M)^{\imath} v^*$. In particular, since this inequality holds for $\imath\rightarrow \infty$, $v^*\leq \omega^*$.
\par
On the other hand, define the operator $\mathscr{M}:\mathcal{V} \rightarrow \mathcal{V}$ for $v\in \mathcal{V}$ by
$
\mathscr{M}v:=\max_{d\in \mathcal{D}}\max_{P_d\in \mathcal{P}_d}
\left\{
Q_d^{-1}r_{(d,P_d)}+Q_d^{-1}R_d v
\right\}.
$
{\color{red}
It is easy to show by a similar argument to the proof of Lemma~\ref{lm1}(a)} that $\mathscr{M}$ is a contraction mapping on $\mathcal{V}$ and the sequence generated by it converges in norm to $v^*$. Moreover, from the definition of $\mathscr{U}^M$ in (\ref{eq23}), $\omega^*=\mathscr{U}^M \omega^* \leq \mathscr{M}^M \omega^*$ holds for $\forall M\in \mathbb{N}$ and especially it is true for $M \rightarrow \infty$. Consequently, $\omega^*\leq v^*$. Therefore, $\omega^*=v^*$.
The proof is completed.


\subsection{Proof of Lemma~\ref{lm5}}
\label{app_E}
{\color{red}
Given a $v$-improving decision rule $d_v\in \mathcal{D}_v\subseteq \mathcal{D}$ and a $d_v$-decreasing transition probability matrix $P_{d_v}^*\in \mathcal{P}_{d_v}^*\subseteq \mathcal{P}_{d_v}$ for $v\in\mathcal{V}$, } we first have $\mathscr{W}^M v= \Phi(d_v,P^{*}_{d_v},v)$
{\color{red}
from (\ref{eqn22}), the definition of $\mathscr{W}^M$ in~(\ref{eq22n}), and the definition of $\Phi(\cdot)$ in~(\ref{eq23}).
}Moreover, from the definition of $\mathscr{U}^M$ in (\ref{eq23}), one have $\mathscr{U}^M u\geq  \Phi(d_v,P^{*}_{d_v},u)$. Then, subtracting $\mathscr{W}^M v$ from $\mathscr{U}^M u$ and utilizing $u\geq v$ lead to
\begin{equation*}
\begin{split}
\mathscr{U}^M u - \mathscr{W}^M v &\geq \Phi(d_v,P^{*}_{d_v},u) -\Phi(d_v,P^{*}_{d_v},v)
\\
&
= ({Q_{d_v}^*}^{-1} R_{d_v}^*)^{M+1} (u-v) \geq 0
,
\end{split}
\end{equation*}
where $(Q_{d_v}^*,R_{d_v}^*)$ is the GS regular splitting of $I-\lambda P_{d_v}^*$.
\par
Furthermore, if $\mathscr{B}u\geq 0$, then from the definition of $\mathscr{W}^M$ in~(\ref{eq22n}) and the fourth equality in~(\ref{eqn22}), one can get
$$
\mathscr{W}^M u = u + \sum_{\varsigma=0}^{M}\left({Q_{d_{u}}^*}^{-1} R_{d_{u}}^*\right)^{\varsigma}( \mathscr{B}u )  \geq u+ \mathscr{B}u
   = \mathscr{Y} u,
$$
where $({Q_{d_{u}}^*}, R_{d_{u}}^*)$ is the GS regular splitting of $I-\lambda P_{d_{u}}^*$ and $P_{d_{u}}^*\in \mathcal{P}_{d_u}^*$. Since for any given $d\in\mathcal{D}$ and $P_d\in \mathcal{P}_d$, $Q_{d}^{-1}r_{(d,P_{d})} +
  Q_{d}^{-1}R_{d} u \geq Q_{d}^{-1}r_{(d,P_{d})} +
  Q_{d}^{-1}R_{d} v$ holds for any regular splitting $(Q_{d},R_d)$ of $I-\lambda P_d$ in view of $u\geq v$, and the maximin operation will not change the direction of the inequality,
  $\mathscr{Y} u\geq \mathscr{Y} v$. Consequently, $\mathscr{W}^M u \geq \mathscr{Y} v$.

\subsection{Proof of Lemma~\ref{lm6}}
\label{app_F}
  Let $u=\mathscr{W}^M v$ and $d_v$ be the $v$-improving decision rule. Then, from Proposition~\ref{prop1},
  there exists a $P'_{d_v}\in \mathcal{P}_{d_v}$ such that
  \begin{equation}\label{F1}
  \mathscr{B}u\geq \mathscr{B}v +({Q'_{d_v}}^{-1}R'_{d_v}-I)(u-v)
  ,
  \end{equation}
   where $(Q'_{d_v},R'_{d_v})$ is the GS regular splitting of $I-\lambda P'_{d_v}$.
  Moreover, let $P_{d_v}^*\in \mathcal{P}_{d_v}^*\subseteq \mathcal{P}_{d_v}$ be the $d_v$-decreasing transition probability matrix. {\color{red}From the definition of $\mathscr{W}^M$ in~(\ref{eq22n}) and the fourth equality in~(\ref{eqn22})}, one can then get
  \begin{equation}\label{F2}
  \begin{split}
  u-v &= \mathscr{W}^M v -v
  =v + \sum_{\varsigma=0}^{M}\left({Q_{d_{v}}^*}^{-1} R_{d_{v}}^*\right)^{\varsigma}\left( \mathscr{B}v \right)
  \\
  &
  -v
  = \sum_{\varsigma=0}^{M}\left({Q_{d_{v}}^*}^{-1} R_{d_{v}}^*\right)^{\varsigma} \left( \mathscr{B}v \right) \geq 0,
  \end{split}
  \end{equation}
  where $(Q_{d_v}^*,R_{d_v}^*)$ is the GS regular splitting of $I-\lambda P_{d_v}^*$. Substituting~(\ref{F2}) into~(\ref{F1}) leads to
  \begin{equation}\label{F3}
  \begin{split}
  &\mathscr{B}u\geq \mathscr{B}v +({Q'_{d_v}}^{-1}R'_{d_v}-I)\sum_{\varsigma=0}^{M}\left({Q_{d_{v}}^*}^{-1} R_{d_{v}}^*\right)^{\varsigma} \left( \mathscr{B}v \right)
  \\
  &
  \geq
  \mathscr{B}v
  + \min_{P\in \mathcal{P}_{d_v}}
  \left(Q^{-1}R - I\right)
  \min_{P\in \mathcal{P}_{d_v}}
  \left\{\sum_{\varsigma=0}^{M}\left(Q^{-1} R\right)^{\varsigma} \left( \mathscr{B}v \right)
  \right\},
  \end{split}
  \end{equation}
where $(Q,R)$ is the GS regular splitting of $I-\lambda P$ for $P\in \mathcal{P}_{d_v}$. Note that for any given $v'\in \mathcal{V}$, both $\left(Q^{-1}R - I\right)v'$ and $\sum_{\varsigma=0}^{M}\left(Q^{-1} R\right)^{\varsigma} v'$ have the same minimum point subject to $P\in \mathcal{P}_{d_v}$. Thus, there exists a $\hat{P}_{d_v}\in \mathcal{P}_{d_v}$ such that
{\color{red}the second term on the right-hand side of~(\ref{F3}) is equal to}
$\left(\hat{Q}_{d_v}^{-1}\hat{R}_{d_v} - I\right)\sum_{\varsigma=0}^{M}\left({\hat{Q}_{d_{v}}}^{-1} \hat{R}_{d_{v}}\right)^{\varsigma} \left( \mathscr{B}v \right)= \left({\hat{Q}_{d_{v}}}^{-1} \hat{R}_{d_{v}}\right)^{M+1} \left( \mathscr{B}v\right) - \mathscr{B}v$,  where $(\hat{Q}_{d_{v}},\hat{R}_{d_{v}})$ is the GS regular splitting of $I-\lambda \hat{P}_{d_v}$. As a result,
$$
\mathscr{B}u \geq \mathscr{B}v + \left({\hat{Q}_{d_{v}}}^{-1} \hat{R}_{d_{v}}\right)^{M+1} \left( \mathscr{B}v\right) - \mathscr{B}v\geq 0.
$$

\subsection{Proof of Theorem~\ref{thm4}}
\label{app_G}
We establish part~(a) by a similar induction argument to Lemma~\ref{lm3}. Note first that $v^0=\tilde{v}^0$ follows from the assumption, and $u^{0}_0(s)$ and $\tilde{u}^{0}_0(s)$ for $\forall s\in \mathcal{S}$ are calculated by the same equations~(\ref{eqn9}),~(\ref{eqA1}), and~(\ref{eqA2}) as in the degenerated case raTVI. Therefore, applying Lemma~\ref{lm3} leads to
\begin{equation}\label{G1}
u^0_{0}-\theta\mathbf{1}\leq \tilde{u}^0_{0} \leq u^0_{0}+\theta\mathbf{1}.
\end{equation}
Moreover, since $\rho_{(s,a)}(v^0)=\rho_{(s,a)}(\tilde{v}^0)$ for $\forall (s,a)\in \mathcal{S}\times\mathcal{A}$ follows from~(\ref{eqn9}) in view of $v^0=\tilde{v}^0$, and for any given $s\in\mathcal{S}$, the approximation error between $\tilde{\rho}_{(s,a)}(\tilde{v}^0)$ and $\rho_{(s,a)}(\tilde{v}^0)$ for $\forall a\in \mathcal{A}$ in~(\ref{eqA1}) does not affect the selection of maximizing actions, $\arg\max_{a\in\mathcal{A}}\tilde{\rho}_{(s,a)}(\tilde{v}^0)
=\arg\max_{a\in\mathcal{A}}\rho_{(s,a)}(v^0)$ for $\forall s\in\mathcal{S}$. Without loss of generality, we select $\tilde{d}_1(s)=d_1(s)$ for $\forall s\in \mathcal{S}$. Then, from~(\ref{eqn9}), one can get that the set of minimizing transition probability distributions are the same for computing $\rho_{(s,a)}(v^0)$ and $\rho_{(s,a)}(\tilde{v}^0)$, $\forall s\in\mathcal{S}$ because $a=\tilde{d}_1(s)=d_1(s)$ and $v^0=\tilde{v}^0$. As such,
{\color{red}
utilizing~(\ref{eqn6}) and the monotonicity of $\rho_{(s,a)}(\cdot)$ shown in Appendix~B, one can get by substituting each term of~(\ref{G1}) that,
}
\begin{equation}\label{G2}
\begin{split}
&\rho_{(s^1,a)}(u^0_{0}-\theta\mathbf{1})=\rho_{(s^1,a)}(u^0_{0})-\lambda\theta\leq \rho_{(s^1,a)}(\tilde{u}^0_{0})
\\
&
\leq \rho_{(s^1,a)}(u^0_{0})+\lambda\theta
=\rho_{(s^1,a)}(u^0_{0}+\theta\mathbf{1}),
\end{split}
\end{equation}
for $a=\tilde{d}_1(s^1)=d_1(s^1)$. In view of $u_1^0(s^1)=\rho_{(s^1,a)}(u^0_{0})$, $a=d_1(s)$ from (\ref{eqA3}) and (\ref{eqA4}) for $\delta=0$, it follows from~(\ref{G2}) that
\begin{equation}\label{G3}
u_1^0(s^1)-\lambda\theta\leq \rho_{(s^1,a)}(\tilde{u}^0_{0})\leq u_1^0(s^1)+\lambda\theta.
\end{equation}
On the other hand, in the inexact case, one can get $\tilde{\rho}_{(s^1,a)}(\tilde{u}^0_{0})-\lambda\delta\leq \rho_{(s^1,a)}(\tilde{u}^0_{0})\leq \tilde{\rho}_{(s^1,a)}(\tilde{u}^0_{0})+\lambda\delta$ for $a=\tilde{d}_1(s)$ from~(\ref{eqA3}), and $\tilde{u}^0_{1}(s^1)=\tilde{\rho}_{(s^1,a)}(\tilde{u}^0_{0})$ from~(\ref{eqA4}). Then, we further have
\begin{equation}\label{G4}
\tilde{u}^0_{1}(s^1)-\lambda\delta \leq \rho_{(s^1,a)}(\tilde{u}^0_{0}) \leq \tilde{u}^0_{1}(s^1)+\lambda\delta.
\end{equation}
{\color{red}
Consequently, combining~(\ref{G3}) with~(\ref{G4}) yields }
$$
u_1^0(s^1)-\theta\leq \tilde{u}_1^0(s^1)\leq u_1^0(s^1)+\theta,
$$
in view of $\lambda(\theta+\delta)=\theta$. Using this inequality and $u_0^0(s^j)-\theta\leq \tilde{u}_0^0(s^j)\leq u_0^0(s^j)+\theta$, $j=2,3,\ldots,m$ from~(\ref{G1}), one can further get
$$
u_1^0(s^2)-\theta\leq \tilde{u}_1^0(s^2)\leq u_1^0(s^2)+\theta,
$$
from~(\ref{eqn6}), (\ref{eqA3}), and (\ref{eqA4}) by a similar argument. By recursion, one can then obtain $u_1^0(s^j)-\theta\leq \tilde{u}_1^0(s^j)\leq u_1^0(s^j)+\theta$, $j=3,4,\ldots,m$ by applying the aforementioned argument for $s^j$. Therefore, $u^0_{1}-\theta\mathbf{1}\leq \tilde{u}^0_{1} \leq u^0_{1}+\theta\mathbf{1}$. Similarly, leveraging this inequality and applying the same argument for $\varsigma=2,3,\ldots,M_0$ by recursion, one can accordingly get $u^0_{\varsigma}-\theta\mathbf{1}\leq \tilde{u}^0_{\varsigma} \leq u^0_{\varsigma}+\theta\mathbf{1}$. Consequently, the induction hypothesis is satisfied for $t=0$.
\par
Suppose that the induction hypothesis is satisfied for $t=\kappa$, and we now show that it also holds for $t=\kappa+1$. First, from the hypothesis for $t=\kappa$, we have
\begin{equation}\label{G5}
u_\varsigma^\kappa-\theta\mathbf{1}\leq \tilde{u}_\varsigma^\kappa
  \leq u_\varsigma^\kappa+\theta\mathbf{1}, \varsigma=0,1,\ldots,M_\kappa.
\end{equation}
Moreover, according to the step~3(f) in Algorithm~1, one can get $v^{\kappa+1}=u_{M_\kappa}^\kappa$ in the exact case and
 $\tilde{v}^{\kappa+1}=\tilde{u}_{M_\kappa}^\kappa$ in the inexact case. Therefore, (\ref{G5}) implies
 $
 v^{\kappa+1}-\theta\mathbf{1} \leq \tilde{v}^{\kappa+1} \leq
 v^{\kappa+1}+\theta\mathbf{1}
 $.
 Applying a similar argument process to the above case for $t=0$, one can first obtain $u^{\kappa+1}_0-\theta\mathbf{1}\leq \tilde{u}^{\kappa+1}_0 \leq u^{\kappa+1}_0+\theta\mathbf{1}$ based on Lemma~\ref{lm3}, and then by recursion, one can further show that $u_\varsigma^{\kappa+1}-\theta\mathbf{1}\leq \tilde{u}_\varsigma^{\kappa+1}
  \leq u_\varsigma^{\kappa+1}+\theta\mathbf{1}$ holds for $\varsigma=1,2,\ldots,M_{\kappa+1}$. As a result, the induction hypothesis holds for $t=\kappa+1$. Thus, part (a) is established.

\par
We next prove part (b). First, when $\delta=0$, one can obtain from the proof of Theorem~\ref{thm3} that there exists a positive integer $N$ for any $\epsilon_1>0$ such that the termination condition $\|u_0^N-v^N\|<(1-\lambda)\epsilon_1/2\lambda$
is satisfied, and meanwhile
\begin{equation}\label{G6}
\begin{split}
\|v^{N+1}-v^*\|
&=\|(\mathscr{T}_{(d^\epsilon,P_{d^\epsilon}^*)})^{M_N+1}v^N-v^*\|
\\
&
\leq \|\mathscr{T}_{(d^\epsilon,P_{d^\epsilon}^*)}v^N-v^*\|
<\epsilon_1/2,
\end{split}
\end{equation}
where the first inequality follows from $v^N\leq \mathscr{Y}v^N=\mathscr{T}_{(d^\epsilon,P_{d^\epsilon}^*)}v^N\leq (\mathscr{T}_{(d^\epsilon,P_{d^\epsilon}^*)})^2v^N\leq \cdots\leq (\mathscr{T}_{(d^\epsilon,P_{d^\epsilon}^*)})^{M_N+1}v^N=v^{N+1}\leq v^*$ by applying Theorem~\ref{thm3} and $v^N\in\mathcal{V}_{\mathscr{B}}$, and the second inequality follows from~(\ref{eq25}).
Moreover, from part (a), we have
$
\|\tilde{u}_0^N-\tilde{v}^N\| \leq \|\tilde{u}_0^N-u_0^N \|+\|u_0^N-v^N\|+\|v^N-\tilde{v}^N\|
\leq \|u_0^N-v^N\|+2\theta <(1-\lambda)\epsilon_1/2\lambda
+2\lambda \delta/(1-\lambda).
$
Given that $\epsilon_1$ is arbitrary, we therefore select a specific $\epsilon_1$ such that
\begin{equation}\label{G7}
\epsilon_1\leq \epsilon -2\lambda(1+\lambda)\delta/(1-\lambda)^2,
\end{equation}
for any $\epsilon>0$. As a result, $\|\tilde{u}_0^N-\tilde{v}^N\|<(1-\lambda)\epsilon/2\lambda
-\delta$. It implies that the raTPI algorithm in the inexact case will also terminate at $t=N$. Suppose that the algorithm terminates with returning a policy $({\tilde{d}}^\epsilon)^\infty$, and denote its corresponding value function under the worst-case transition probability matrix by $\tilde{v}^\epsilon$. From~(\ref{eqn23}), one can see that $\tilde{v}^\epsilon$ can be regard as the value $\tilde{v}^{N+1}$ generated by Algorithm~1 at $t=N+1$ by selecting $M_{N}\rightarrow \infty$ and $\delta=0$ when $t=N$. Therefore, from part~(a), (\ref{G6}), and~(\ref{G7}), one can get $
\|v^\epsilon - v^*\|=\|\tilde{v}^{N+1} - v^*\|\leq \|\tilde{v}^{N+1} -v^{N+1}\|+\|v^{N+1}-v^*\|< \theta+ \epsilon_1/2< \epsilon
$. The proof is completed.

\bibliographystyle{ieeetr}
\bibliography{ref}

\end{document}